\documentclass[a4paper,10pt]{article}
\usepackage{anysize}
\marginsize{3.0cm}{3.0cm}{2.5cm}{1.5cm}
\usepackage{amsmath,amssymb,amsthm,ulem}
\usepackage{latexsym,amsfonts,amsmath}
\usepackage{graphicx,psfrag}
\usepackage{epstopdf}
\usepackage{pstricks}
\usepackage{amssymb}
\usepackage{amsthm}

\def\NN{\mathbb{N}}

\def\RR{\mathbb{R}}

\def\liminf{\mathop{\underline{\lim}}}

\newcommand{\defi}{\stackrel{\triangle}{=}}

\newtheorem{condition}{Condition}{\bfseries}{\itshape}
{\bfseries}{\itshape}

\newtheorem{theorem}{Theorem}{\bfseries}{\itshape}

\newtheorem{corollary}{Corollary}{\bfseries}{\itshape}

\newtheorem{proposition}{Proposition}{\bfseries}{\itshape}

{\bfseries}{\itshape}

{\bfseries}{\itshape}

\newtheorem{remark}{Remark}{\bfseries}{\itshape}

\newtheorem{definition}{Definition}{\bfseries}{\itshape}

\begin{document}
\author{
Alexey Piunovskiy \\
Department of Mathematical Sciences, University of Liverpool, L69 7ZL, UK.\\ \texttt{piunov@liv.ac.uk}\\ \ \\
Alexander Plakhov\\
Center for R\&{}D in Mathematics and Applications, Department of Mathematics,\\ University of Aveiro, 3810-193, Portugal\\
and Institute for Information Transmission Problems, Moscow, 127051, Russia\\ \texttt{plakhov@ua.pt}\\ \ \\
Delfim F. M. Torres\\
Center for R\&{}D in Mathematics and Applications, Department of Mathematics,\\ University of Aveiro, 3810-193, Portugal\\  \texttt{delfim@ua.pt}\\ \ \\
Yi Zhang\\
Department of Mathematical Sciences, University of Liverpool, L69 7ZL, UK.\\
\texttt{Yi.Zhang@liv.ac.uk}
}
\title{Optimal {Impulse}  Control of Dynamical Systems}
\date{}

\maketitle

\begin{abstract}
Using the tools of the Markov Decision Processes, we justify the dynamic programming approach to the optimal impulse control of deterministic dynamical systems. We prove the equivalence of the integral and differential forms of the optimality equation.  The theory is illustrated by an example from mathematical epidemiology. The developed methods can be also useful for the study of piecewise deterministic Markov processes.
\end{abstract}
\begin{tabbing}
\small \hspace*{\parindent}  \= {\bf Keywords:}
Dynamical System, Impulse Control, Total Cost, Discounted Cost, \\ Randomized Strategy, {Piecewise Deterministic Markov Process}\\
\> {\bf AMS 2000 subject classification:} \= Primary 49N25; Secondary 90C40.
\end{tabbing}

\section{Introduction}
{Impulse control of various dynamical systems attracts attention of many researchers: \cite{b14,b1,b2,b3,b4,b5,b6,b7,b10,b11,b15,b13,b8,b16,b12,b9}, to mention the most relevant and the most recent works. The underlying system can be described in terms of ordinary \cite{b14,b1,b2,b10,b11,b13,b8} or stochastic \cite{b15} differential equations; that was an abstract Markov process in \cite{b16}. In \cite{b3,b4,b5,b7,b12,b9}, along with the given deterministic drift, there are spontaneous (or natural) Markov jumps of the state. Such models are called Piecewise Deterministic Markov Processes (PDMP); the drift is usually described by a fixed flow. On the other hand, if there is no drift and the trajectories are piecewise constant, the model is called  Continuous-Time Markov Decision Process (MDP) \cite{b7}. By the gradual control we mean that only the local characteristics of the underlying process are under control. In case of PDMP, it means that the deterministic drift and the rate of the spontaneous/natural jumps, as well as the post-jump distribution are under control. But the {\it impulse} control means the following: at particular discrete time moments, the decision maker decides to intervene by instantaneously moving the process to some new point in the state space; that new point may be also random. Then, restarting at this new point, the process runs until the next intervention and so on. Sometimes, such control is called `singular control' \cite{b15}. The goal is to minimize the total (expected) accumulated cost which may be discounted \cite{b1,b3,b4,b5,b6,b7,b10,b15,b16,b12,b9} or not \cite{b14,b1,b2,b10,b11,b8,b9}. The case of long-run average cost was also studied in e.g. \cite{b9}.

The most popular method of attack to such problems is Dynamic Programming \cite{b1,b3,b4,b5,b6,b7,b15,b16,b12,b9}. In \cite{b10,b11,b8}, versions of the Pontryagin Maximum Principle is used. In \cite{b2}, the impulse control is firstly reformulated as the linear program on measures: impulses correspond to the singular, Dirac components. After that, the numerical approximate scheme is developed in the form of Linear Matrix Inequalities.

Impulse control theory is widely applied to different real-life problems: epidemiology \cite{b14,b13}, Internet congestion control \cite{b1}, reliability \cite{b5}, economic and finance \cite{b10,b15,b12}, moving objects \cite{b10}, medicine \cite{b11},  genetics and ecology \cite{b8} etc.

The distinguishing features of the current work are as follows.
\begin{itemize}
\item We consider the purely deterministic positive model with the total cost. As is known and explained in the text, the discounted model is a special case, as well as the absorbing model.
\item The imposed conditions partially overlap with those introduced in other articles. Generally speaking, our conditions are weaker than the assumptions introduced in the cited literature.
\item For the model under study, we demonstrate the new method to obtain the optimality equation in the integral form and to develop the corresponding successive approximations. This method is based on the well known tools for Discrete-Time MDP.
\item Under mild conditions, we prove the equivalence of the optimality (Bellman) equation in the integral and differential form. The analytical proof is new. Moreover, as mentioned in Conclusion, this proof remains valid also for the more general case of PDMP. Note also that the differential form is slightly different from what appeared in other works.
\item We present the solution to the optimal impulse control of an epidemic model, which is of its own interest.
\end{itemize}

The paper is organized in the following way. After describing the problem statement, we demonstrate the MDP approach and provide the integral optimality equation in  Section \ref{sec3}. In Section \ref{sec4}, we prove the equivalence of the integral and differential forms of the optimality equation. The impulse control of SIR epidemic is developed in Section \ref{sec:SIR}. In Conclusion, we briefly describe the ways for generalizing our results to PDMP.}

The following notations are frequently used throughout this paper. $\NN=\{1,2,\ldots\}$ is the set of natural numbers;
$\delta_{x}(\cdot)$
is the Dirac measure concentrated at $x$, we call such distributions degenerate; $I\{\cdot\}$ is the indicator function.
${\cal{B}}(E)$ is the Borel $\sigma$-algebra of the Borel space $E$, ${\cal P}(E)$ is the Borel space of probability measures on $E$. (It is always clear which $\sigma$-algebra is fixed in $E$.) The Borel $\sigma$-algebra ${\cal B}({\cal P}(E))$ comes from the weak convergence of measures, after we fix a proper topology in $E$. $\mathbb{R}_+\defi (0,+\infty)$, $\mathbb{R}_+^0\defi[0,+\infty)$, $\bar\RR^0_+\defi[0,+\infty]$; in $\RR_+$ and $\RR_+^0$, we consider the Borel $\sigma$-algebra, and $Leb$ is the Lebesgue measure.
The abbreviation $w.r.t.$ (resp. $a.s.$) stands for ``with respect to'' (resp. ``almost surely''); for $b\in[-\infty,+\infty]$, $b^+\defi\max\{b,0\}$ and $b^-\defi\min\{b,0\}$. $\inf\emptyset\defi +\infty$. Measures introduced in the current article can take infinite values. { If $b=\infty$ then the integrals $\displaystyle \int_{(a,b]}f(u)du$ are taken over the open interval $(a,\infty)$.}

\section{Problem Statement}\label{sec2}

We will deal with a control model defined through the following elements.
\begin{itemize}
\item $\bf X$ is the state space, a Borel subset of a complete separable metric space with metric $\rho_X$ and the Borel $\sigma$-algebra.
\item $\phi(\cdot\,, \cdot):~{\bf X}\times\RR^0_+\to{\bf X}$ is the flow possessing the semigroup property $\phi(x,t+s)=\phi(\phi(x,s),t)$ for all $x\in{\bf X}$ and $(t,s)\in(\RR^0_+)^2$; $\phi(x,0)=x$ for all $x\in{\bf X}$. Between the impulses, the state changes according to the flow.
\item $\bf A$ is the action space, again a Borel subset of a complete separable metric space with metric $\rho_A$ and the Borel $\sigma$-algebra.
\item $l(\cdot\,, \cdot):~{\bf X}\times{\bf A}\to{\bf X}$ is the mapping describing the new state after the corresponding action/impulse is applied.
\item $C^g(\cdot):~{\bf X}\to\RR$ is the (gradual) cost rate.
\item $C^I(\cdot\,, \cdot):~{\bf X}\times{\bf A}\to\RR$ is the cost associated with the actions/impulses applied in the corresponding states.
\end{itemize}
All the mappings $\phi(\cdot),l(\cdot),C^g(\cdot)$ and $C^I(\cdot)$ are assumed to be measurable.

Let ${\bf X}_\Delta={\bf X}\cup\{\Delta\}$, where $\Delta$ is an isolated artificial point describing the case that the controlled process is over and no future costs will appear. The dynamics (trajectory) of the system can be represented as one of the following sequences
\begin{eqnarray}
&&x_0\to (\theta_1,a_1)\to x_1\to (\theta_2,a_2)\to \ldots;~~~~ \theta_i<+\infty \mbox{ for all } i\in\NN, \nonumber\\
\mbox{or}&&\label{e1}\\
&&x_0\to (\theta_1,a_1)\to\ldots\to x_n\to (+\infty,a_{n+1})\to \Delta\to (\theta_{n+2},a_{n+2}) \to \Delta \to \ldots,\nonumber
\end{eqnarray}
where $x_0\in{\bf X}$ is the initial state of the controlled process and $\theta_i<+\infty$ for all $i=1,2,\ldots, n$. For the state $x_{i-1}\in{\bf X}$, $i\in\NN$, the pair $(\theta_i,a_i)\in\bar\RR^0_+\times{\bf A}$ is the control at the step $i$: after $\theta_i$ time units, the impulsive action $a_i$ will be applied leading to the new state
\begin{equation}\label{e1p}
x_i=\left\{\begin{array}{ll}
l(\phi(x_{i-1},\theta_i),a_i), & \mbox{ if } \theta_i<+\infty; \\
\Delta, & \mbox{ if } \theta_i=+\infty.
\end{array}\right.
\end{equation}
The state $\Delta$ will appear forever, after it appeared for the first time, i.e., it is absorbing.

After each impulsive action, if $\theta_1,\theta_2,\ldots,\theta_{i-1}<+\infty$, the decision maker has in hand complete information about the history, that is, the sequence
$$x_0, (\theta_1,a_1), x_1,\ldots,  (\theta_{i-1},a_{i-1}),x_{i-1}.$$
The next control $(\theta_i,a_i)$ is based on this information and we allow the pair $(\theta_i,a_i)$ to be randomized.

The cost on the coming interval of the length $\theta_i$ equals
\begin{equation}\label{e1pprim}
\int_{(0,\theta_i] } C^g(\phi(x_{i-1},u))du+I\{\theta_i<+\infty\} C^I(\phi(x_{i-1},\theta_i),a_i),
\end{equation}
the last term being absent if $\theta_i=+\infty$.
If the cost functions $C^g(\cdot)$ and $C^I(\cdot)$ can take positive and negative values, then one can calculate separately the expressions in (\ref{e1pprim}) for the positive and negative parts of the costs and accept the convention $+\infty-\infty\defi +\infty$.
The next state $x_i$ is given by formula (\ref{e1p}).

In the space of all the trajectories (\ref{e1})
\begin{eqnarray*}
\Omega&=&\bigcup_{n=1}^\infty[{\bf X}\times((\RR^0_+\times{\bf A})\times{\bf X})^n
\times(\{+\infty\}\times{\bf A})\times\{\Delta\}\times ((\bar\RR^0_+\times{\bf A})\times\{\Delta\})^\infty]\\
&& \bigcup [{\bf X}\times((\RR^0_+\times{\bf A})\times{\bf X})^\infty],
\end{eqnarray*}
we fix the natural $\sigma$-algebra $\cal F$.
Finite sequences
$$h_i=(x_0, (\theta_1,a_1), x_1, (\theta_2,a_2),\ldots,x_i)$$
will be called (finite) histories; $i=0,1,2,\ldots$, and the space of all such histories will be denoted as ${\bf H}_i$. Capital letters $X_i,T_i,\Theta_i, A_i$ and $H_i$ denote the corresponding functions of $\omega\in\Omega$, i.e., random elements.

\begin{definition}\label{d1}
A control strategy $\pi=\{\pi_i\}_{i=1}^\infty$ is a sequence of stochastic kernels $\pi_i$ on $\bar\RR^0_+\times{\bf A}$ given ${\bf H}_{i-1}$. A control strategy is called stationary deterministic and denoted as $(\varphi_\theta,\varphi_a)$, if, for all $i=1,2,\ldots$, $\pi_i(d\theta\times da|h_{i-1})=\delta_{\varphi_\theta(x_{i-1})}(d\theta) \delta_{\varphi_a(x_{i-1})}(da)$, where $\varphi_\theta:{\bf X}_\Delta\to\bar\RR^0_+$ and $\varphi_a:{\bf X}_\Delta\to{\bf A}$ are measurable mappings.
\end{definition}

If the initial state $x_0\in{\bf X}$ and a strategy $\pi$ are fixed, then there is a unique probability measure $P^\pi_{x_0}(\cdot)$ on $\Omega$ satisfying the following conditions:
$$P^\pi_{x_0}(X_0\in\Gamma_X)=\delta_{x_0}(\Gamma_X)\mbox{ for } \Gamma_X\in{\cal B}({\bf X}_\Delta);$$
for all $i\in\NN$, $\Gamma\in{\cal B}(\bar\RR^0_+\times{\bf A})$, $\Gamma_X\in{\cal B}({\bf X}_\Delta)$,
\begin{eqnarray*}
P^\pi_{x_0}((\Theta_i,A_i)\in\Gamma|H_{i-1})&=&\pi_i(\Gamma |H_{i-1});\\
P^\pi_{x_0}(X_i\in\Gamma_X|H_{i-1},(\Theta_i,A_i))
&=& \left\{\begin{array}{ll}
\delta_{l(\phi(X_{i-1},\Theta_i),A_i)}(\Gamma_X), & \mbox{ if } X_{i-1}\in{\bf X},~\Theta_i<+\infty; \\
\delta_\Delta(\Gamma_X) & \mbox{ otherwise } \end{array} \right.
\end{eqnarray*}
For details, see the Ionescu Tulcea Theorem \cite[Prop.7.28]{Bertsekas:1978}. The mathematical expectation w.r.t. $P^\pi_{x_0}$ is denoted as $E^\pi_{x_0}$.

The optimal control problem under study looks as follows.
\begin{eqnarray}
\mbox{Minimize w.r.t. } \pi &&\nonumber\\
{\cal V}(x_0,\pi) &=&  E^\pi_{x_0}\left[\sum_{i=1}^\infty I\{X_{i-1}\ne\Delta\} \left\{ \int_{(0,\Theta_{i}]}  C^g(\phi(X_{i-1},u)) du\right.\right.\label{PZZeqn02}\\
&&+ I\{\Theta_i<+\infty\} \left.\left.\vphantom{\sum_{i=1}^\infty} C^I(\phi(X_{i-1},\Theta_i),A_i)\right\}\right].\nonumber
\end{eqnarray}

\begin{definition}\label{d2}
A control strategy $\pi^*$ is called uniformly optimal if, for all $x\in{\bf X}_\Delta$, ${\cal V}(x,\pi^*)={\cal V}^*(x)\defi\inf_{\pi}{\cal V}(x,\pi)$.
\end{definition}

\section{MDP Approach}\label{sec3}
In this section, we establish the optimality results for problem (\ref{PZZeqn02}) by referring to the known ones for its induced total undiscounted Markov Decision Process (MDP).

The MDP under study is given by the state and action spaces ${\bf X}_\Delta$ and $\bar\RR^0_+\times{\bf A}$, transition kernel
$$Q(dy|x,(\theta,a))\defi\left\{\begin{array}{ll}
\delta_{l(\phi(x,\theta),a)}(dy), & \mbox{ if } x\ne\Delta,~\theta\ne+\infty;\\ \delta_\Delta(dy) & \mbox{ otherwise}, \end{array}\right. $$
and cost
$$\tilde C(x,(\theta,a))\defi I\{x\ne\Delta\} \left\{ \int_{(0,\theta]}  C^g(\phi(x,u)) du
+ I\{\theta<+\infty\} C^I(\phi(x,\theta),a)\right\}.$$
Clearly, actions of the form $(+\infty,a)$ can be treated as stopping the control process with the terminal cost $\displaystyle \int_{(0,+\infty)} C^g(\phi(x,u))du$. In this framework, we denote as {\it "stop"} the strategy which immediately chooses $\theta_1=+\infty$: ${\cal V}(x_0,"stop")= \int_{(0,+\infty)} C^g(\phi(x,u))du$.
The artificial state $\Delta$ means that MDP is stopped without any future cost. { If,  for all $x\in{\bf X}$,  $\inf_\pi{\cal V}(x_0,\pi)<+\infty$, then the optimal control problem (\ref{PZZeqn02}) is not degenerate. This assumption holds if $\displaystyle\int_{(0,\infty)} C^g(\phi(x,u))du<+\infty$ or, e.g., in the following cases.}

\begin{itemize}
\item { Absorbing case: there is a specific measurable  "cemetery" subset ${\bf Y}\subset {\bf X}$ such that for all $y\in{\bf Y}$ $C^g(y)\equiv 0$,
 for all $u\in\RR^0_+$ $\phi(y,u)\in{\bf Y}$ and, for each $x\in{\bf X}\setminus{\bf Y}$,\\
 -- either $\inf_{u\in\RR^0_+} \{u:~\phi(x,u)\in{\bf Y}\}\defi t^*(x)<+\infty$, the function $t^*(x)$ is measurable, and $\displaystyle \int_{(0,t^*(x)]} C^g(\phi(x,u))du<+\infty$,\\
 -- or there is an action $\hat a\in{\bf A}$ such that $l(x,\hat a)\in{\bf Y}$. }
\item Discounted case: the state space $\bf X$ has the form ${\bf X}={\bf Y}\times\RR^0_+$, where $\bf Y$ is a Borel subset of a complete separable metric space, the component $s_0$ of the initial state $x_0=(y_0,s_0)$ is zero, the flow $\phi(\cdot)$ satisfies $\phi((y,s),t)=(\phi_{\bf Y}(y,t),s+t)$, where $\phi_{\bf Y}(\cdot)$ is a flow in $\bf Y$, and
    \begin{eqnarray*}
    C^g((y,s))&=&e^{-\alpha s}C^g_{\bf Y}(y);\\
    C^I((y,s),a)&=&e^{-\alpha s} C^I_{\bf Y}(y,a);\\
    l((y,s),a)&=&(l_{\bf Y}(y,a),s).
    \end{eqnarray*}
    Here the functions and mapping $C^g_{\bf Y}(\cdot)$, $C^I_{\bf Y}(\cdot)$ and $l_{\bf Y}(\cdot)$ are for the component $y\in{\bf Y}$ only, and, for all $y\in{\bf Y}$, $\displaystyle \int_{(0,+\infty)} e^{-\alpha u} C^g_{\bf Y}(\phi_{\bf Y}(y,u))du<+\infty$. $\alpha>0$ is the discount factor, component $s$ of the state $x=(y,s)$ plays the role of time, and in principle one can consider the non-homogeneous model with the functions and mappings $C^g_{\bf Y}(\cdot)$, $C^I_{\bf Y}(\cdot)$ and $l_{\bf Y}(\cdot)$ depending also on the component $s$.
\item Generalized discounting: the model is as in the previous item, but
$$C^g((y,s))=h(s)C^g_{\bf Y}(y);~~~~~C^I((y,s),a)=h(s) C^I_{\bf Y}(y,a),$$
and the measurable function $h(\cdot)$ is such that $\forall y\in{\bf Y}$ $\displaystyle \int_{(0,+\infty)} h(s) C^g_{\bf Y}(\phi_{\bf Y}(y,u))du<+\infty$.
\end{itemize}

Throughout this section, the following condition is satisfied.
\begin{condition}\label{asa}
$C^g$ and $C^I$ are $\mathbb{R}_+^0$-valued, that is, we consider the so called positive model with the total expected cost.
\end{condition}

This condition means that we deal with a positive model. In this case, the value function ${\cal V}^*(x_0)\defi \inf_{\pi}{\cal V}(x_0,\pi)$ is the minimal $\bar{\mathbb{R}}_+^0$-valued lower semianalytic solution to the following optimality (Bellman) equation:
 \begin{eqnarray}
 V(\Delta)&=&0,\nonumber\\
 V(x)&=&\inf_{(\theta,a)\in\bar{\mathbb{R}}_+^0\times\textbf{A}}\left\{\tilde C(x,(\theta,a))+\int_{{\bf X}_\Delta} V(y) Q(dy|x,(\theta,a))\right\}\label{PZZeqn01}\\
 &=&\inf_{(\theta,a)\in\bar{\mathbb{R}}_+^0\times\textbf{A}}\left\{\int_{(0,\theta]}  C^g(\phi(x,u)) du
+ I\{\theta<+\infty\}\left( C^I(\phi(x,\theta),a)+V(l(\phi(x,\theta),a))\right)\right\}\nonumber \\
&& ~~~~    \forall~x\in\textbf{X}.\nonumber
 \end{eqnarray}
(See \cite[Cor.9.4.1,Prop.9.8, and Prop.9.10]{Bertsekas:1978}.)

{Recall that our model is the special case of PDMP when the spontaneous (natural) jumps intensity $\lambda$ equals zero. In case of discounted cost, corresponding versions of equation (\ref{PZZeqn01}) appeared in the works \cite{b3,b4,b5,b7,b12} on PDMP.}

Note that the case of a simultaneous sequence of impulses, when $\Theta_i=\Theta_{i+1}=\ldots=0$, is not excluded. In such cases, the total cost (\ref{PZZeqn02}) is calculated over a finite time horizon, up to the accumulation of impulses.

Section \ref{sec4} is devoted to the study of equation (\ref{PZZeqn01}),  without requiring that Condition \ref{asa} is satisfied.

In the framework of stopping MDP, the decision to stop (here that means $\theta=+\infty$, and all the actions of the form $(+\infty,a)$ can be merged to one point) is usually considered as an isolated point of the action space $\bar\RR^0_+\times{\bf A}$. But in this  case the remainder (real) action space $\RR^0_+\times{\bf A}$ would be not compact. To avoid this inconvenience, we accept the following conditions.

\begin{condition}\label{PPZcondition01}
\begin{itemize}
\item[(a)] The space $\bf A$ is compact, and $+\infty$ is the one-point compactification of the positive real line $\RR^0_+$, so that the action space $\bar\RR^0_+\times{\bf A}$ in the MDP is compact.
\item[(b)] The mapping $(x,a)\in \textbf{X}\times \textbf{A}\rightarrow l(x,a)$ is continuous.
\item[(c)] The mapping $(x,\theta)\in \textbf{X}\times \mathbb{R}_+^0\rightarrow \phi(x,\theta)$ is continuous.
\item[(d)] The function $(x,a)\in \textbf{X}\times \textbf{A}\rightarrow C^I(x, a)$ is lower semicontinuous.
\item[(e)] The function $x\in \textbf{X}\rightarrow C^g(x)$ is lower semicontinuous.
\end{itemize}
\end{condition}

Still under these conditions, the model is not semicontinuous because, if $x\ne\Delta$, $\theta_n\in\RR^0_+$ and $\theta_n\to+\infty$ then the transition probabilities $Q(dy|x,(\theta_n,a))$ do not converge to $\delta_\Delta(dy)$. Nevertheless, the usual dynamic programming approach is fruitful.

\begin{theorem}\label{t1}
Suppose Conditions \ref{asa} and \ref{PPZcondition01} are satisfied. Then the following assertions hold.
\begin{itemize}
\item[(1)] The minimal $\bar\RR^0_+$-valued solution $V(x)$ to equation (\ref{PZZeqn01}) is lower semicontinuous, unique, and can be constructed by successive approximations starting from $V_0(x)\equiv 0$, $x\in{\bf X}_\Delta$:
    \begin{eqnarray*}
    V_{n+1}(\Delta)&=& V_n(\Delta);\\
    V_{n+1}(x)&=& \inf_{(\theta,a)\in\bar{\mathbb{R}}_+^0\times\textbf{A}}\left\{\int_{(0,\theta]}  C^g(\phi(x,u)) du
+ I\{\theta<+\infty\}\left( C^I(\phi(x,\theta),a)\right.\right.\\
&&\left.\vphantom{\int_t }\left.+V_n(l(\phi(x,\theta),a))\right)\right\}~~~~~    \forall~x\in\textbf{X}.
    \end{eqnarray*}
The sequence $\{V_n\}_{n=1}^\infty$ increases point-wise and $V(x)=\lim_{n\to\infty} V_n(x)={\cal V}^*(x)$.
\item[(2)] There exist measurable mappings $\varphi^*_\theta:{\bf X}\to\bar\RR^0_+$ and $\varphi^*_a:{\bf X}\to{\bf A}$ such that, for all $x\in{\bf X}$,
\begin{eqnarray}\label{e2}
{\cal V}^*(x)&=&\int_{(0,\varphi^*_\theta(x)]}  C^g(\phi(x,u)) du\\
&&+ I\{\varphi^*_\theta(x)<+\infty\}\left( C^I(\phi(x,\varphi^*_\theta(x)),\varphi^*_a(x))
+{\cal V}^*(l(\phi(x,\varphi^*_\theta(x)),\varphi^*_a(x)))\right).
\nonumber
\end{eqnarray}
\item[(3)] A stationary deterministic strategy $(\varphi^*_\theta,\varphi^*_a)$ is uniformly optimal if and only if it satisfies equality (\ref{e2}).
\end{itemize}
\end{theorem}

\underline{Proof.} It is sufficient to consider only $x\in{\bf X}$, as $\Delta$ is the isolated point of ${\bf X}_\Delta$. Suppose $W(\cdot)$ is a lower semicontinuous $\bar\RR^0_+$-valued function on $\bf X$ and show that function on ${\bf X}\times\bar\RR^0_+\times{\bf A}$
\begin{equation}\label{e3}
\int_{(0,\theta]}  C^g(\phi(x,u)) du
+ I\{\theta<+\infty\}\left( C^I(\phi(x,\theta),a)+W(l(\phi(x,\theta),a))\right)
\end{equation}
is lower semicontinuous and $\bar\RR^0_+$-valued.

Firstly, let us show that the non-negative function $\displaystyle \int_{(0,\theta]} C^g(\phi(x,u))du$ is lower semicontinuous on ${\bf X}\times\bar\RR^0_+$. By a well known result of Baire, see e.g., \cite[Lemma 7.14]{Bertsekas:1978}, there exists an increasing sequence of bounded $\mathbb{R}_+^0$-valued continuous functions, say $\{c_m(\cdot)\}_{m=1}^\infty$, on $\textbf{X}$ such that $c_m(x)\uparrow C^g(x)$ for each $x\in \textbf{X}_.$
For every $m=1,2,\dots,$ function $(x,\theta)\in {\bf X}\times\bar\RR^0_+\rightarrow \int_{(0,\theta]} c_m(\phi(x,u))e^{-\frac{u}{m}} du $ is bounded continuous. By the monotone convergence theorem and using the result of Baire again, we see that function
\begin{eqnarray*}
(x,\theta)\in {\bf X}\times\bar\RR^0_+ \rightarrow \int_{(0,\theta]} C^g(\phi(x,u))du &=&\int_{(0,\theta]} \lim_{m\to\infty} c_m(\phi(x,u))e^{-\frac{u}{m}} du\\
&=&\lim_{m\to\infty} \int_{(0,\theta]} c_m(\phi(x,u))e^{-\frac{u}{m}} du
\end{eqnarray*}
is lower semicontinuous.

Secondly, let us show that the non-negative function
$$F(x,\theta,a)\defi I\{\theta<+\infty\}\left( C^I(\phi(x,\theta),a)+W(l(\phi(x,\theta),a))\right)$$
is also lower semicontinuous on ${\bf X}\times\bar\RR^0_+\times{\bf A}$. Suppose $(x_n,\theta_n,a_n)\to (x,\theta\ne+\infty,a)$ as $n\to\infty$. Since the flow $\phi(\cdot)$ and the mapping $l(\cdot)$ are continuous, we deduce that the sequences $y_n\defi \phi(x_n,\theta_n)$ and $l_n\defi l(\phi(x_n,\theta_n),a_n)$ converge to $y\defi\phi(x,\theta)$ and $l\defi l(\phi(x,\theta),a)$ correspondingly. Therefore
$$\liminf_{n\to\infty} F(x_n,\theta_n,a_n)\ge\liminf_{n\to\infty} C^I(y_n,a_n)+\liminf_{n\to\infty}W(l_n,a_n)\ge C^I(y,a)+W(l,a)= F(x,\theta,a)$$ because the both functions $C^I(\cdot)$ and $W(\cdot)$ are lower semicontinuous. (See \cite[Lemma 7.13]{Bertsekas:1978}). In case $(x_n,\theta_n,a_n)\to (x,+\infty,a)$ as $n\to\infty$, it is obvious that $\liminf_{n\to\infty} F(x_n,\theta_n,a_n)\ge 0=F(x,+\infty,a)$.

Therefore, function (\ref{e3}) is lower semicontinuous and obviously $\bar\RR^0_+$-valued.

(1) Clearly, $V_1(x)\ge V_0(x)$, so that the sequence $\{V_n\}_{n=1}^\infty$ increases point-wise and hence converges to some $\bar\RR^0_+$-valued function $V(\cdot)$. The non-negative function $V_0(\cdot)$ is lower semicontinuous and, if $V_n(\cdot)$ is a non-negative lower semicontinuous function then so is  function $V_{n+1}(\cdot)$ by Proposition 7.32 \cite{Bertsekas:1978}. Therefore, function $V(\cdot)$ is lower semicontinuous by the mentioned above Baire result.

For every $n=0,1,2,\ldots$, function (\ref{e3}), with $W(\cdot)$ being replaced with $V_n(\cdot)$, is lower semicontinuous. Therefore, the set
\begin{eqnarray*}
U_n(x,\lambda)&\defi &\{(\theta,a)\in\bar\RR^0_+\times{\bf A}:~\int_{(0,\theta]} C^g(\phi(x,u))du+I\{\theta<+\infty\} (C^I(\phi(x,\theta),a)\\
&&+V_n(l(\phi(x,\theta),a)))\le\lambda\}
\end{eqnarray*}
is closed and hence compact for all $x\in{\bf X}$, $\lambda\in\RR$. (See Condition \ref{PPZcondition01}(a).) By Proposition 9.17 \cite{Bertsekas:1978}, $V(x)={\cal V}^*(x)$.

(2) The value function ${\cal V}^*(\cdot)$ satisfies equation (\ref{PZZeqn01}) and is lower semicontinuous. By Proposition 7.33 \cite{Bertsekas:1978}, there exists a measurable mapping $\varphi^*:~{\bf X}\to\bar\RR^0_+\times{\bf A}$ which provides the infimum in (\ref{PZZeqn01}) for all $x\in{\bf X}$. Assertion (2) follows.

(3) This assertion follows directly from Proposition 9.12 \cite{Bertsekas:1978}.

~\hfill $\Box$

\begin{corollary}\label{cor3}
Suppose Conditions \ref{asa} and \ref{PPZcondition01} are satisfied and $(\varphi^*_\theta,\varphi^*_a)$ is a uniformly optimal stationary deterministic strategy. Then the following assertions hold true.
\begin{itemize}
\item[(1)] For every $x\in{\bf X}$ with $\varphi^*_\theta(x)<+\infty$, for $y\defi\phi(x,\varphi^*_\theta(x))$, equality
$$\inf_{a\in{\bf A}}\{C^I(y,a)+V(l(y,a))\}=C^I(y,\varphi^*_a(x))+V(l(y,\varphi^*_a(x)))=V(y)$$
is valid.
\item[(2)] For every $x\in{\bf X}$ with $\varphi^*_\theta(x)>0$, for each $t\in[0,\varphi^*_\theta(x))$, for $y\defi\phi(x,t)$, equality
\begin{eqnarray}
V(y) &=&
\int_{(0,\varphi^*_\theta(x)-t]} C^g(\phi(y,u))du+I\{\varphi^*_\theta(x)<+\infty\}[C^I(\phi(y,\varphi^*_\theta(x)-t),\varphi^*_a(x)) \nonumber\\
&& +V(l(\phi(y,\varphi^*_\theta(x)-t),\varphi^*_a(x)))] \label{e26}
\end{eqnarray}
\end{itemize}
is valid and hence
\begin{equation}\label{e28}
V(x)=\int_{(0,t]} C^g(\phi(x,u))du+V(\phi(x,t)).
\end{equation}
\end{corollary}

\underline{Proof.} (1) The left equality is obvious.

The case when $V(y)$ is bigger than the expression of the left is excluded. (Consider $\theta=0$ for $y$ in (\ref{PZZeqn01}).) If $V(y)$ is smaller then, the pair
$$(\tilde\varphi_\theta(x)\defi \varphi^*_\theta(x)+\varphi^*_\theta(y),~\varphi^*_a(y))$$
gives rise to the  expression
\begin{eqnarray*}
&& \int_{(0,\varphi^*_\theta(x)]} C^g(\phi(x,u))du+\int_{(\varphi^*_\theta(x),\tilde\varphi_\theta(x)]} C^g(\phi(x,u))du\\
&&+I\{\tilde\varphi_\theta(x)<+\infty\}[C^I(\phi(x,\tilde\varphi_\theta(x)),\varphi^*_a(y))+ V(l(\phi(x,\tilde\varphi_\theta(x)),\varphi^*_a(y)))]\\
&=& \int_{(0,\varphi^*_\theta(x)]} C^g(\phi(x,u))du+ V(y),
\end{eqnarray*}
that is, $(\tilde\varphi_\theta(x),\varphi^*_a(y)=\varphi^*_a(\phi(x,\varphi^*_\theta(x))))$
provides the smaller value for
\begin{equation}\label{e27}
\int_{(0,\theta]} C^g(\phi(x,u))du+I\{\theta<+\infty\} [C^I(\phi(x,\theta),a)+V(l(\phi(x,\theta),a))]
\end{equation}
than the pair $(\varphi^*_\theta,\varphi^*_a)$, which contradicts the definition of $(\varphi^*_\theta,\varphi^*_a)$.

(2) If $\varphi^*_\theta(x)<+\infty$, then $\varphi^*_a(x)$ provides the infimum to
\begin{eqnarray*}
&C^I(\phi(y,\varphi^*_\theta(x)-t),a)+V(l(\varphi(y,\varphi^*_\theta(x)-t),a))
=C^I(\phi(x,\varphi^*_\theta(x)),a)+V(l(\varphi(x,\varphi^*_\theta(x)),a)).
\end{eqnarray*}
$V(y)$ cannot be bigger than the expression on the left in (\ref{e26}). (Consider $\theta=\varphi^*_\theta(x)-t$ for $y$ in (\ref{PZZeqn01}).) If $V(y)$ is smaller then, like previously, the pair
$$(\tilde\varphi_\theta(x)=t+\varphi^*_\theta(y),~\varphi^*_a(y))$$
provides the smaller value for (\ref{e27}) than the pair $(\varphi^*_\theta,\varphi^*_a)$, which contradicts the definition of $(\varphi^*_\theta,\varphi^*_a)$.

Equality (\ref{e28}) follows from (\ref{e26}) because
\begin{eqnarray*}
V(x)&=&\int_{(0,t]} C^g(\phi(x,u))du+\int_{(0,\varphi^*_\theta(x)-t]} C^g(\phi(y,u))du\\
&&+I\{\varphi^*_\theta(x)<+\infty\}[C^I(\phi(y,\varphi^*_\theta(x)-t),\varphi^*_a(x)) +V(l(\phi(y,\varphi^*_\theta(x)-t),\varphi^*_a(x)))].
\end{eqnarray*}

~\hfill $\Box$

\begin{remark}\label{rem11}
If Conditions \ref{asa} and \ref{PPZcondition01} are satisfied then, for each $x\in{\bf X}$, in case the pair $(\hat\theta<+\infty,\varphi^*_a(x))$ provides the infimum in equation (\ref{PZZeqn01}), equality
\begin{equation}\label{EqR1}
C^I(\phi(x,\hat\theta),\varphi^*_a(x))+V(l(\phi(x,\hat\theta),\varphi^*_a(x)))=V(\phi(x,\hat\theta))
\end{equation}
is valid. The proof coincides with the proof of Item (1) of Corollary \ref{cor3}.
\end{remark}

Below, it will be assumed that the function ${\cal V}^*(\cdot)$ is finite-valued. This requirement is obviously satisfied for positive models if, for each $x\in{\bf X}$, there exists a control strategy $\pi$ such that ${\cal V}(x,\pi)<+\infty$. The latter assumption follows from the following condition.

\begin{condition}\label{con2} For all $x \in \mathbf{X}$ the composite function $C^g(\phi(x,t))$ is Lebesgue integrable on $\mathbb{R}_+^0$. This means that the integral
$$
\int_{(0,{+\infty})} C^g(\phi(x,t))\, dt
$$
exists and is finite.
\end{condition}

In what follows,
we accept the following convention.
We say that a function $g : \bf X \to \mathbb{R}$ satisfies a certain property (is continuous, absolutely continuous, measurable, Lebesgue integrable, etc.) along the flow $\phi$, if for all $x \in \bf X$ the composite function $t \mapsto g(\phi(x,t))$ from $\mathbb{R}_+^0$ to $\mathbb{R}$ satisfies this property. In view of this convention, Condition \ref{con2} asserts that the function $C^g$ is Lebesgue integrable along the flow.

The following proposition states that for each $x\in{\bf X}$ the set of values $\theta$ providing the infimum in (\ref{PZZeqn01})  is closed in $\bar\RR^0_+$, and hence, contains its minimal value denoted as $\theta^*(x)$.

Under Conditions \ref{asa} and \ref{PPZcondition01}, for the minimal non-negative solution $V$ to equation (\ref{PZZeqn01}),
we introduce the function $G:~{\bf X}\times\bar\RR^0_+\times{\bf A}\to\bar\RR^0_+$  by
$$
G(x,\theta,a)\defi \int_{(0,\theta]} C^g(\phi(x,u))du+I\{\theta<+\infty\}\left[C^I(\phi(x,\theta),a)+V(l(\phi(x,\theta),a))\right]
$$
and the sets  $\Theta(x)$ by
\begin{equation}\label{e30}
\qquad \Theta(x) \defi \left\{\theta\in\bar\RR^0_+:~\inf_{a\in{\bf A}} G(x,\theta,a) =V(x)\right\}.
\end{equation}
For a fixed $x\in{\bf X}$, the set $\Theta(x)$ contains all time moments $\theta$ such that the pair $(\theta,\hat a)$ provides the infimum in (\ref{PZZeqn01}). Here, for $\theta\in\Theta(x)\cap\RR^0_+$, $\hat a\in{\bf A}$ provides the infimum in (\ref{e30}); such $\hat a$ exist because the function $G$ is lower semicontinuous in $a$, if Conditions \ref{asa} and \ref{PPZcondition01} are satisfied.

\begin{proposition}\label{prop0}
Suppose Conditions \ref{asa} and \ref{PPZcondition01} are satisfied. Then for all $x \in \mathbf{X}$ the set $\Theta(x)$ is non-empty and closed in $\bar\RR^0_+$, and therefore contains the value $\theta^*(x) \defi \inf \Theta(x)$, the mapping $\theta^*(\cdot)$ being measurable.
\end{proposition}

The proofs of all propositions are postponed to the Appendix.

Recall that, under Conditions \ref{asa} and \ref{PPZcondition01}, there is a measurable mapping $\hat\varphi^*_a:~{\bf X}\times\RR^0_+\to{\bf A}$ providing
\begin{equation} \label{e25}
\inf_{a\in{\bf A}}\left\{C^I(\phi(x,\theta),a)+V(l(\phi(x,\theta),a))\right\}.
\end{equation}
(See Proposition 7.33 of \cite{Bertsekas:1978}.)
The mapping
$\theta^*(x) = \inf\Theta(x)$ from Proposition \ref{prop0} is measurable, so that the pair $(\varphi^*_\theta,\varphi^*_a)$ with $\varphi^*_\theta(x)\defi \theta^*(x)$ and
$\varphi^*_a(x)\defi\hat\varphi^*_a(x,\theta^*(x))$ satisfies Item (2) of Theorem \ref{t1}.

Further, we will need the following two conditions strengthening Conditions  \ref{con2} and \ref{asa} respectively.
\begin{condition}\label{as3} There is $K \in \mathbb{R}_+$ such that for all $x\in{\bf X}$, $\displaystyle \int_{(0,+\infty)} |C^g(\phi(x,u))|du\le K$.
\end{condition}

\begin{condition}\label{asb} The function $C^g$ is $\mathbb{R}_+^0$-valued and $C^I\ge \delta>0$.
\end{condition}

Condition \ref{asb}  guarantees that, for reasonable strategies $\pi$, $\Theta_i$ is finite only a finite number of times ($P^\pi_{x_0}$-a.s.): otherwise ${\cal V}(x_0,\pi)=+\infty$.

Under Conditions \ref{as3} and \ref{asb}, starting from any initial state $x_0\in{\bf X}$, for the optimal strategy the MDP must be stopped at a finite time moment
$$T_{stop}=\min\{i:~\Theta_{i+1}=+\infty\}$$
and, for any reasonable strategy $\pi$,
$$E^\pi_{x_0} [T_{stop}]\le \frac{K}{\delta}:$$
otherwise, the total cost
$${\cal V}(x_0,\pi)=E^\pi_{x_0}\left[\sum_{i=0}^\infty \tilde C(X_i,(\Theta_i,A_i))\right]>K$$
is bigger than that coming from stopping the MDP immediately:
$${\cal V}(x_0,"stop")=\int_{(0,+\infty)} C^g(\phi(x,u))du\le K.$$
If we restrict ourselves to such control strategies, then we are in the framework of absorbing MDP \cite[\S9.6]{las99}, and the following proposition can be proved similarly to Theorem 9.6.10(c) \cite{las99}.

\begin{proposition}\label{prop3}
Suppose Conditions \ref{PPZcondition01}, \ref{as3} and \ref{asb} are satisfied. Then the Bellman equation (\ref{PZZeqn01}) has a unique bounded lower semicontinuous solution.
\end{proposition}

{ \begin{remark} \label{rem9}
According to the proof of Proposition \ref{prop3} (see Appendix),
the Bellman equation (\ref{PZZeqn01}) has a unique bounded lower semicontinuous solution also in the case when Conditions \ref{asa} and \ref{PPZcondition01} are satisfied, the function ${\cal V}^*$ is bounded and, for all strategies $\pi$, for all $x_0\in{\bf X}$, $\lim_{i\to\infty} E^\pi_{x_0}[V(X_i)]=0$, where $V$ is a bounded lower semicontinuous function satisfying equation (\ref{PZZeqn01}). {In this case $V = \mathcal{V}^*$.}
\end{remark} }


\section{Differential Form of the Optimality Equation}\label{sec4}

In this section, we establish the equivalence of the integral and differential formulations of the optimality equation using minimal assumptions about the system. We do not assume any structure of the sets $\mathbf{X}$ and $\mathbf{A}$ and of the maps $C^I : \bf X \times \mathbf{A} \to \mathbb{R}$ and $l : \bf X \times \mathbf{A} \to \bf X$; we only require Condition \ref{con2} to be satisfied.  Firstly, we justify the differential form of the optimality equation for the general model studied in Section \ref{sec3}. After that, we briefly  discuss the discounted model.

\subsection{Total Cost Model}\label{ssec41}

\begin{definition}\label{def3}
A point $x$ is said to be a {\it singular point} of the flow $\phi$, if it is not an intermediate point of a trajectory, that is, the equation $x = \phi(\tilde x, s)$ has no solutions for all $\tilde x \in \bf X$ and $s > 0$. Note that if the flow possesses the {\it group property}, i.e., $\phi(x,t+s)=\phi(\phi(x,s),t)$ for all $x\in{\bf X}$ and $(t,s)\in \RR^2$, then there are no singular points.
\end{definition}

Let $V : \mathbf{X} \to \mathbb{R}$ be a certain function. We denote
$$
\mathcal{F}_+^V(x) \defi \lim_{t \to 0^+} \Big[ \frac{V(\phi(x,t)) - {V(x)}}{t} + \frac{1}{t} \int_{(0,t]} C^g(\phi(x,u))\, du \Big],
$$
provided that the limit in the right hand side exists.

Further, if $x$ is a nonsingular point of the flow, we define the number set
\begin{eqnarray*}
\underline{\cal{F}}_-^V(x)& \defi & \left\{ \underline{\lim}_{t \to 0^+} \Big[ \frac{V(x) - V(\phi(\tilde{x}, s-t))}{t} \right.
\\
&&\left. + \frac{1}{t} \int_{(-t,0]} C^g(\phi(\tilde{x}, s+u))\, du \Big] :~~(\tilde x,s)\in{\bf X}\times\RR_+, ~\phi(\tilde x,s)=x\right\} \subset \mathbb{R} \cup \{ \pm\infty \}.
\end{eqnarray*}
With some abuse of notation, if $\underline{\cal{F}}_-^V(x)$ is a singleton (e.g. if the flow possesses the group property), then we identify it with its element.
If, otherwise, $x$ is a singular point, then we set $\underline{\cal{F}}_-^V(x) = \emptyset $.

\begin{remark}\label{rem2}
{If ${\bf X}'\subset \mathbf{X}\subset \RR^d$ is} a smooth open manifold, the flow is given by the differential equation $\dot x = f(x)$, satisfying the standard conditions on the existence of a unique local solution in ${\bf X}'$ (for positive and negative $t$), for each initial condition from ${\bf X}'$, and $C^g(x)$ is continuous along the flow in ${\bf X}'$ and $V(x)$ is continuously differentiable along the flow in ${\bf X}'$, then $\underline{\cal{F}}_-(x)$ is a singleton for all $x\in{\bf X}'$ and
$$
\mathcal{F}_+^V(x) = \underline{\cal{F}}_-^V(x) = C^g(x) + \nabla V(x) \cdot f(x).
$$\end{remark}

Consider the optimality equation (\ref{PZZeqn01}) on $\bf X$, that is,
the following integral equation:
\begin{equation}\label{IntEq}
V(x) = \inf_{\theta \in \bar{R}_+^0} \Big\{ \int_{(0,\theta]} C^g(\phi(x,u))\, du + I\{ \theta < +\infty \}\, \inf_{a \in \mathbf{A}} \{ C^I(\phi(x,\theta), a) + V(l(\phi(x,\theta), a)) \} \Big\},~~~x\in{\bf X}.
\end{equation}
Everywhere further, we assume that the function $V(x)$ is finite-valued. For example, under Conditions \ref{asa}, \ref{PPZcondition01} and \ref{con2} the value function ${\cal V}^*(x)=V(x)$, studied in Section \ref{sec3}, is finite and satisfies the equation (\ref{IntEq}) by Theorem \ref{t1}.

\begin{condition} \label{con41} For each $x$ the former infimum in the right hand side of (\ref{IntEq}) is attained on a nonempty set $\Theta(x) \subset \bar{\mathbb{R}}_+^0$, and $\Theta(x)$ contains its infimum.
\end{condition}

We emphasise that Condition \ref{con41} is satisfied under Conditions \ref{asa} and \ref{PPZcondition01} if $\bf X$ is a Borel space: see Proposition \ref{prop0}.

We also consider the so called Bellman equation in the differential form:
$$~$$

for all $x\in \mathbf{X}$,
\begin{equation}\label{DifEq1}
\begin{split}
 \hspace*{-20mm}  \text{either \hspace*{4mm} (a) } \hspace*{18mm}  &\mathcal{F}_+^V(x) = 0
 \\
 \text{and} \quad &\inf_{a\in\mathbf{A}} \big[ C^I(x, a) + V(l(x, a)) - V(x) \big] > 0,
\\
 \hspace*{-20mm} \text{or  \hspace*{10mm} (b) }  \hspace*{18mm}  &\underline{\cal{F}}_-^V(x) \subset \bar{\mathbb{R}}_+^0
\\
\text{and} \quad &\inf_{a\in\mathbf{A}} \big[ C^I(x, a) + V(l(x, a)) - V(x) \big] = 0.
\end{split}
\end{equation}

\begin{remark}\label{rem5}
In the case (a) it is assumed that the right {limit} {\rm exists} and equals 0. If it does not exist then the case (b) should take place.
\end{remark}

{Recall that our model is the special case of PDMP when the spontaneous (natural) jumps intensity $\lambda$ equals zero. In case of discounted cost, corresponding versions of equation (\ref{DifEq1}) appeared in the works \cite{b3,b5,b7,b12,b9} on PDMP, see also the paper \cite{b1} on the purely deterministic system. In \cite{b9}, the undiscounted case was also investigated.
More about connection of the current work with existing results at the end of Subsection \ref{ssec42}. Here, we only emphasize that the differential form in the shape of inclusion $\underline{\cal{F}}_-^V(x) \subset \bar{\mathbb{R}}_+^0$ did not appear in the cited literature. Remember, $\underline{\cal{F}}_-^V(x)$ is a singleton in case the flow possesses the group property.}

Define the set
\begin{equation}\label{e29}
\mathcal{L} \defi \{ x \in {\bf X} : \inf_{a\in\mathbf{A}} \big[ C^I(x, a) + V(l(x, a)) - V(x) \big] = 0 \}.
\end{equation}
If $V$ is  a solution to equation (\ref{PZZeqn01}), the set $\cal L$
can be understood as the set of the states, where actions/impulses must be applied. Suppose Conditions \ref{asa} and \ref{PPZcondition01} are satisfied and function $V$ is the minimal $\bar\RR^0_+$-valued solution to equation (\ref{PZZeqn01}). For each $x\in{\bf X}$, if $\hat\theta\in\Theta(x)\cap\RR^0_+$ (see (\ref{e30})), then $\hat\theta$ provides the infimum to $\inf_{a\in{\bf A}} G(x,\theta,a)$ and hence the pair $(\hat\theta,\varphi^*_a(x))$, where $C^I(\phi(x,\hat\theta),\varphi^*_a(x))+V(l(\phi(x,\hat\theta),\varphi^*_a(x)))=V(\phi(x,\hat\theta))$, provides the infimum in (\ref{PZZeqn01}). According to Remark \ref{rem11}, $\phi(x,\hat\theta)\in{\cal L}$, that is, $\Theta(x)\cap\RR^0_+\subset\{t\in\RR^0_+:~\phi(x,t)\in{\cal L}\}$. Usually, this inclusion is strict, and $\Theta(x)\cap\RR^0_+$ is a singleton coinciding with the infimum of $\{t\in\RR^0_+:~\phi(x,t)\in{\cal L}\}$.

We consider the following conditions on the function $V(x)$ satisfying equation (\ref{DifEq1}).
\vspace{2mm}

\begin{condition}\label{con4} For all $x$ the set $\{ t \in \mathbb{R}_+^0 : \phi(x,t) \in \mathcal{L} \}$ is either empty, or contains its infimum.
\end{condition}

\begin{condition}\label{con5}  The function $V(x)$ is right lower semicontinuous and left upper semicontinuous along the flow. That is, first, for all $x$ we have
$$
\underline{\lim}_{t \to 0^+} V(\phi(x,t)) \ge V(x).
$$
Second, for all $x$ and all $(\tilde x,s)\in{\bf X}\times\RR_+$ such that $\phi(\tilde x, s) = x$ we have
$$
\overline{\lim}_{t \to 0^+} V(\phi(\tilde x, s-t)) \le V(x).
$$
(If $x$ is singular, the inequality is satisfied by default.)
\end{condition}

\begin{condition}\label{con6}  If, for some $x \in \mathbf{X}$ and $s > 0$ and for all $0 \le t < s$ the states $\phi(x,t)$ are not in $\mathcal {L}$,
then $\lim_{t \to s^-} V(\phi(x,t)) = V(\phi(x,s))$. In other words, if the relative interior points of the flow trajectory between $x$ and $\phi(x,s)$ are not contained in $\mathcal{L}$ then $V(x)$ is left continuous along the flow at $\phi(x,s)$.
\end{condition}

{In the following theorem, we establish the equivalence of the integral and differential forms of the optimality equation. To the best of our knowledge, such an analytical proof never appeared in the existing literature. Note that we do not assume that the flow possesses the group property. }

\begin{theorem}\label{th1} Suppose Condition \ref{con2} is satisfied.
Let the  function $V : \mathbf{X} \to \mathbb{R}$ be measurable along the flow, and additionally, the integral $\displaystyle\int_{(0,+\infty)} V(\phi(x,t))\, dt$ be finite for all $x\in{\bf X}$. Then the following statements are equivalent.
\begin{itemize}
\item[(1)] $V(x)$ satisfies equation (\ref{IntEq}) and Condition \ref{con41};
\item[(2)] $V(x)$ satisfies equation (\ref{DifEq1}) and Conditions \ref{con4}-\ref{con6}.
\end{itemize}
\end{theorem}

\underline{Proof.} Below, we use the notation
$$\mathcal{I}V(x) \defi  \inf_{a \in \mathbf{A}} \{ C^I(x, a) + V(l(x, a)) \}.$$

{\bf 1.} Suppose that assertion (1) is valid and prove assertion (2).

Let $x\in{\bf X}$ be fixed. For any $t > 0$ we can write down
\begin{equation}\label{a1}
V(x) \le \inf_{\theta \in [t,\,+\infty]} \Big\{ \int_{(0,\theta]} C^g(\phi(x,u))\, du + I\{ \theta < +\infty \}\, \mathcal{I}V(\phi(x, \theta)) \Big\}.
\end{equation}
Changing the variables $u - t = v, \ \theta - t = s$, using the semigroup property $\phi(x, v+t) = \phi(\phi(x,t), v)$ and denoting for brevity $x' = \phi(x,t)$, we get
\begin{equation}\label{a2}
\begin{split}
V(x) \le \int_{(0,t]} C^g(\phi(x,u))\, du + \inf_{s\in\bar{\mathbb{R}}_+^0} \Big\{ \int_{(0,s]} C^g(\phi(x', v))\, dv
\\
+ I\{ s < +\infty \}\, \mathcal{I}V(\phi(x', s)) \Big\},
\end{split}
\end{equation}
and thus,
\begin{equation}\label{a3}
V(x) \le \int_{(0,t]} C^g(\phi(x,u))\, du + V(\phi(x,t)).
\end{equation}
One can easily see that the integral in the right hand side goes to 0 as $t \to 0^+$, and, taking the lower limit of the both parts in this inequality, one obtains
$$
V(x) \le \underline{\lim}_{t \to 0^+} V(\phi(x,t)).
$$
That is, $V$ is right lower semicontinuous along the flow.

If $x$ is not singular, let $x = \phi(\tilde x, s)$ with $s > 0$. For all $t \in [0,\, s]$ we have
$$
V(\phi(\tilde x, s-t)) \le \inf_{\theta \in [t,\,+\infty]} \Big\{ \int_{(0,\theta]} C^g(\phi(\tilde x, s-t+u))\, du
+ I\{ \theta < +\infty \}\, \mathcal{I}V(\phi(\tilde x, s-t+\theta)) \Big\}.
$$
Changing the variables $u - t = v$ and $\theta - t = \tau$, we get
$$
V(\phi(\tilde x, s-t)) \le \int_{(-t,0]} C^g(\phi(\tilde x, s+v))\, dv
$$
$$
+\inf_{\tau\in\bar{\mathbb{R}}_+^0} \Big\{ \int_{(0,\tau]} C^g(\phi(\tilde x, s+v))\, dv
+ I\{ \tau < +\infty \}\, \mathcal{I}V(\phi(\tilde x, s+\tau)) \Big\}.
$$
Taking into account that $\phi(\tilde x, s) = x$ and using (\ref{IntEq}) we obtain
\begin{equation}\label{a4}
V(\phi(\tilde x, s-t)) \le \int_{(-t,0]} C^g(\phi(\tilde x, s+v))\, dv + V(x).
\end{equation}
Now taking the upper limit of the both parts in this relation and using that the integral in the right hand side goes to 0 as $t \to 0^+$, one gets
$$
\overline{\lim}_{t \to 0^+} V(\phi(\tilde x, s-t)) \le V(x).
$$
That is, $V$ is left upper semicontinuous along the flow, which means that Condition \ref{con5} is satisfied.

Recall that $\Theta(x)$ is the (nonempty) set of values $\theta$ minimizing (\ref{IntEq}). Let $\theta^*(x) = \inf \Theta(x)$. By Condition \ref{con41} we have $\theta^*(x) \in \Theta(x)$. Consider two cases: $\theta^*(x) > 0$ and $\theta^*(x) = 0$, and prove equation (\ref{DifEq1}).
 \vspace{2mm}

($\alpha$) $\theta^*(x) > 0$.

Take a finite $t \in [0,\, \theta^*(x)]$. Then (\ref{IntEq}) remains valid if the infimum is taken over $\theta \in [t,\, +\infty]$; as a consequence, we have the equality instead of "$\le\,$" in relations (\ref{a1}) and (\ref{a2}). Moreover, the infimum in (\ref{a2}) is attained at $s\in\Theta(x) - t$. Hence
\begin{equation}\label{formula1m}
V(x) = \int_{(0,t]} C^g(\phi(x,u))\, du + V(\phi(x,t))
\end{equation}
and $\Theta(\phi(x,t)) = \Theta(x) - t$, and therefore,
\begin{equation}\label{formula2}
\theta^*(\phi(x,t)) = \theta^*(x) - t.
\end{equation}

 It follows from (\ref{formula1m}) that
 $$
 \frac{ V(\phi(x,t)) - V(x)}{t} + \frac{1}{t} \int_{(0,t]} C^g(\phi(x,u))\, du = 0,
 $$
and therefore,
$$
\mathcal{F}_+^V(x) = 0.
$$

Now using that the infimum in (\ref{IntEq}) is not attained at $\theta = 0$, we have $V(x) < \mathcal{I}V(x)$, and therefore,
$$
\inf_{a\in\mathbf{A}} \big[ C^I(x, a) + V(l(x, a)) - V(x) \big] > 0.
$$
Thus, relations (\ref{DifEq1}a) are valid, and therefore $x \not\in \mathcal{L}$.
\vspace{1mm}

($\beta$) $\theta^*(x) = 0$.

If $x$ is a singular point of the flow then the first relation in (\ref{DifEq1}b) is transformed into the valid formula $\emptyset \subset \bar{\mathbb{R}}_+^0$.

Suppose $x$ is a nonsingular point, that is, $x = \phi(\tilde x, s)$ with $s > 0$.
Rewrite inequality (\ref{a4}) as follows
$$
\frac{V(\phi(\tilde x, s-t)) -  V(x)}{t} \le \frac{1}{t}  \int_{(-t,0]} C^g(\phi(\tilde x, s+v))\, dv;
$$
hence
$$
\underline{\lim}_{t \to 0^+} \Big[ \frac{ V(x) - V(\phi(\tilde x, s-t))}{t} + \frac{1}{t} \int_{(-t,0]} C^g(\phi(\tilde x, s+v))\, dv \Big] \ge 0.
$$
It follows that
$$
\underline{\cal{F}}_-^V(x) \subset \bar{\mathbb{R}}_+^0.
$$

Further, since $0 \in \Theta(x)$ and therefore the infimum in (\ref{IntEq}) is attained at $\theta = 0$, we have $V(x) = \mathcal{I}V(x) = \inf_{a \in \mathbf{A}} \{ C^I(x, a) + V(l(x, a))$. It follows that
$$
\inf_{a\in\mathbf{A}} \big[ C^I(x, a) + V(l(x, a)) - V(x) \big] = 0.
$$
Thus, relations (\ref{DifEq1}b) are valid, and therefore $x \in \mathcal{L}$.
\vspace{2mm}

Let us check Condition \ref{con6}. Suppose $s > 0$ and $\phi(x,t) \not\in \mathcal{L}$ for all $0 \le t < s$. This means that $\theta^*(\phi(x,t)) > 0$, and using (\ref{formula2}) we conclude that $\theta^*(x)>t$ for all $t\in[0,s)$; hence $\theta^*(x) \ge s$. Substituting { $s$ for $t$} in formula (\ref{formula1m}) we obtain $V(x) = \int_{(0,s]} C^g(\phi(x,u))\, du + V(\phi(x,s))$. Subtracting (\ref{formula1m}) from this formula, one obtains
$$
0 = \int_{(t,s]} C^g(\phi(x,u))\, du + V(\phi(x,s)) - V(\phi(x,t)),
$$
and therefore,
$$
V(\phi(x,t)) =  \int_{(t,s]} C^g(\phi(x,u))\, du+V(\phi(x,s)).
$$
It follows that $\lim_{t \to s^-} V(\phi(x,t)) = V(\phi(x,s))$, and so, Condition \ref{con6} is satisfied.

It remains to check Condition \ref{con4}.
From ($\alpha$) and ($\beta$) we conclude that if $\theta^*(x) > 0$ then $x \not\in \mathcal{L}$, and if $\theta^*(x) = 0$ then $x \in \mathcal{L}$. By formula (\ref{formula2}), if $0 \le t < \theta^*(x)$ then $\theta^*(\phi(x,t)) = \theta^*(x) - t > 0$ and therefore $\phi(x,t) \not\in \mathcal{L}$, and  if $t = \theta^*(x)$ then $\theta^*(\phi(x,t)) = 0$ and so, $\phi(x,t) \in \mathcal{L}$. It follows that if $\theta^*(x) = +\infty$ then the set
$$\{ t \subset {\mathbb{R}}_+^0: \phi(x,t) \in \mathcal{L} \}$$
is empty, and if $\theta^*(x) < +\infty$ then $\theta^*(x)$ is contained in this set and is its infimum. Thus, Condition \ref{con4} is satisfied.
\vspace{2mm}

{\bf 2}. Suppose that assertion (2) is valid and prove assertion (1). In the proof below we use the following statements, which will be proved in Appendix.

\begin{proposition}\label{prop2}
Let $t \in \bar{\mathbb{R}}_+^0$. If $h$ is defined on $[0,\, t] \cap {\mathbb{R}}_+^0$ and for all $s \in (0,\, t)$
$$\mbox{ either }~~\underline{h}'_-(s)\defi \liminf_{t\to 0^+}\frac{h(s)-h(s-t)}{t} \ge 0~~~ \mbox{ or }~~ \underline{h}'_+(s)\defi \liminf_{t\to 0_+}\frac{h(s+t)-h(s)}{t} \ge 0,$$
and additionally, $h$ is right lower semicontinuous on $[0,\, t)$ and left upper semicontinuous on $(0,\, t] \cap {\mathbb{R}}_+^0$, then $h$ is monotone nondecreasing.
\end{proposition}

\begin{proposition}\label{prop1}
Let $t \in \bar{\mathbb{R}}_+^0$. If $h$ is defined on $[0,\, t] \cap {\mathbb{R}}_+^0$ and is left continuous on $(0,\, t] \cap {\mathbb{R}}_+^0$, and $h'_+(s) = 0$ for $s \in [0,\, t)$ then $h$ is constant.
\end{proposition}

Fix arbitrary $x \in \mathbf{X}$ and $t \in \bar{\mathbb{R}}_+^0$ and define the function $h$ by
$$
h(s) = V(\phi(x,s)) -  \int_{(s,t]} C^g(\phi(x, v))\, dv - I\{ t < +\infty \}\, {\mathcal{I}} V(\phi(x,t)),
$$
with $s \in [0,\, t] \cap {\mathbb{R}}_+^0$.

First, we show that $h$ is (monotone) nondecreasing.
For $s\in[0,t)$, consider $\frac{h(s+\tau)-h(s)}{\tau}$ under $\tau\in(0,t-s)$. This ratio equals
$$
\frac{V(\phi(x, s+\tau)) - V(\phi(x, s))}{\tau} +  \frac{1}{\tau} \int_{(s,s+\tau]} C^g(\phi(x,v))\, dv.
$$
Denoting $x' = \phi(x,s)$ and using the semigroup property of the flow and that the integral in the right hand side goes to zero as $\tau \to 0^+$, we get
\begin{equation}\label{2a}
h'_+(s) = \lim_{\tau \to 0^+} \Big\{ \frac{V(\phi(x', \tau)) - V(x')}{\tau} + \frac{1}{\tau} \int_{(0,\tau]} C^g(\phi(x', u))\, du \Big\}={\cal F}^V_+(\phi(x,s)),
\end{equation}
if $h'_+(s)$ and ${\cal F}^V_+(\phi(x,s))$ exist. We emphasize that ${\cal F}^V_+(\phi(x,s))$ and $h'_+(s)$ exist (or do not exist) simultaneously.

In a similar way one calculates the lower left derivative of $h$ and concludes that
\begin{equation}\label{2b}
\underline{h}'_-(s) \in \underline{\cal{F}}_-^V(\phi(x,s))
\end{equation}
 for all $s\in(0,t]\cap\RR$ and also for $s=0$, provided $x$ is not a singular point.

According to equation (\ref{DifEq1}), for $x \not\in \mathcal{L}$ the derivative $\mathcal{F}_+^V(x)$ exists and equals zero, and for $x \in \mathcal{L}$ we have $\underline{\cal{F}}_-^V(x) \subset \bar{\mathbb{R}}_+^0.$ By (\ref{2a}) and (\ref{2b}) we conclude that if $\phi(x, s) \not\in \mathcal{L}$ then $h'_+(s)$ exists and equals zero, and,  if $\phi(x, s) \in \mathcal{L}$ then $\underline{h}'_-(s) \ge 0$. Taking into account Condition \ref{con5}, we conclude that the function $h$ is also right lower { semicontinuous in $[0,\, t)$} and left upper semicontinuous { in $(0,\, t] \cap \mathbb{R}_+^0$}. Therefore $h$ satisfies all conditions of Proposition \ref{prop2} and hence is nondecreasing.

If $t < +\infty$, we have
\begin{equation}\label{2-2}
h(t) = V(\phi(x,t)) - {\mathcal{I}} V(\phi(x,t)),
\end{equation}
and by virtue of (\ref{DifEq1}), $h(t) \le 0$. Thus,
$$
h(0) = V(x) -  \int_{(0,t]} C^g(\phi(x, v))\, dv - {\mathcal{I}} V(\phi(x,t)) \le h(t) \le 0.
$$

If, otherwise, $t = +\infty$, we have for all $s\in\RR^0_+$
$$
h(0) \le h(s) = V(\phi(x,s)) - \int_{(s,+\infty)} C^g(\phi(x, v))\, dv.
$$
According to Condition \ref{con2}, the integral in the right hand side of this relation is finite and goes to 0 as $s \to +\infty$. The function $s \mapsto V(\phi(x,s))$ has a finite Lebesgue integral on ${\mathbb{R}}_+^0$. Therefore { $h(s) \le 0$ for all $s\in\RR^0_+$ and in particular,
$$
h(0) = V(x) -  \int_{(0,+\infty)} C^g(\phi(x,u))\, du \le 0.
$$}

This proves that
\begin{equation}\label{eq1}
V(x) \le \inf_{t\in \bar{\mathbb{R}}_+^0} \left\{ \int_{(0,t]\cap\RR} C^g(\phi(x,u))\, du + I\{ t < +\infty\}\, {\mathcal{I}} V(\phi(x,t)) \right\}.
\end{equation}
\vspace{1mm}

Now for each $x\in{\bf X}$ take the value
$$
t^*(x) := \inf \{ t \in {\mathbb{R}}_+^0 : \phi(x,t) \in \mathcal{L} \} \in \bar{\mathbb{R}}_+^0
$$
and set $t = t^*(x)$ in the definition of the function $h$. We intend to show that in this case $h(s)\equiv 0$.

For $0 \le s < t^*(x)$ we have $\phi(x,s) \not\in \mathcal{L}$ and therefore $h'_+(s) = 0$.
According to Condition \ref{con6}, the function $s \mapsto V(\phi(x, s))$, and therefore also the function $h$, are left continuous for finite $0 < s \le t^*(x)$. Hence by Proposition \ref{prop1}\, $h$ is constant on $[0,\, t^*(x)] \cap {\mathbb{R}}_+^0$.

Let $t^*(x)$ be finite. Condition \ref{con4} states that $\phi(x, t^*(x)) \in \mathcal{L}$, and therefore,
$$
h(t^*(x)) = V(\phi(x, t^*(x))) - {\mathcal{I}} V(\phi(x, t^*(x))) = 0.
$$
It follows that $h(s) = 0$ for all $0 \le s \le t^*(x)$, and in particular,
$$
h(0) = V(x) - \int_{(0,t^*(x)]} C^g(\phi(x,u))\, du - {\mathcal{I}} V(\phi(x,t^*(x))) = 0.
$$

If, otherwise, $t^*(x) = +\infty$, the constant function $h(s)$ is the sum of a function that has a finite Lebesgue integral on ${\mathbb{R}}_+^0$ and a function going to 0 as $s \to +\infty$; therefore it is the null function, and again,
$$
h(0) = V(x) -  \int_{(0,+\infty)} C^g(\phi(x,u))\, du = 0.
$$
This proves that
$$
V(x) =  \int_{(0,t^*(x)]} C^g(\phi(x,u))\, du + I\{ t^*(x) < +\infty \}\, {\mathcal{I}} V(\phi(x,t^*(x)))
$$
\begin{equation}\label{eq2}
\ge \inf_{t\in \bar{\mathbb{R}}_+^0} \left\{ \int_{(0,t]} C^g(\phi(x,u))\, du + I\{ t < +\infty \}\, {\mathcal{I}} V(\phi(x,t)) \right\}.
\end{equation}
As a consequence of (\ref{eq1}) and (\ref{eq2}) we obtain that equation (\ref{IntEq}) is true and $t^*(x) \in \Theta(x)$.

It remains to show that Condition \ref{con41} is satisfied. Set an arbitrary $0 \le t < t^*(x)$ in the definition of $h$. Since $\phi(x,t) \not\in \mathcal{L}$, by (\ref{2-2}) and (\ref{DifEq1}a) we have
$$
h(t) =  V(\phi(x,t)) - {\mathcal{I}} V(\phi(x,t)) < 0,
$$
and taking into account that $h$ is nondecreasing, we conclude that
$$
h(0) = V(x) -  \int_{(0,t]} C^g(\phi(x,u))\, du - {\mathcal{I}} V(\phi(x,t)) \le h(t) < 0.
$$
As a result we have
$$
V(x) < \int_{(0,t]} C^g(\phi(x,u))\, du + {\mathcal{I}} V(\phi(x,t));
$$
that is, $t \not\in \Theta(x)$ for $0 \le t < t^*(x)$. This implies that $t^*(x) = \inf\Theta(x)$, and so, Condition \ref{con41} is satisfied.
\hfill $\Box$

\begin{remark}\label{rem8}
Obviously, for all $x\in{\bf X}$, $t\in\RR_+$, ${\cal V}^*(\phi(x,t))\le \int_{(t,+\infty)} C^g(\phi(x,u))du$. Therefore, if, under Condition \ref{asa}, the following stronger version of Condition \ref{con2}:
$$\forall x\in{\bf X}~~~~~\int_{(0,+\infty)}\left(\int_{(t,+\infty)} C^g(\phi(x,u))du\right) dt<+\infty$$
is satisfied, then the integral $\int_{(0,+\infty)} {\cal V}^*(\phi(x,t))dt$ is finite,  provided the function ${\cal V}^*:~{\bf X}\to\RR$ is Lebesgue-measurable along the flow. By Theorem \ref{t1}, under Conditions \ref{asa} and \ref{PPZcondition01} the latter requirement is satisfied and ${\cal V}^*=V$ is the minimal non-negative solution to equation (\ref{PZZeqn01}).
\end{remark}

\subsection{Discounted Model}\label{ssec42}

Note that for the validity of Theorem \ref{th1}, only Condition \ref{con2} is needed. In the case of the discounted model described in Section \ref{sec3}, it takes the following form.

The function $C^g_{\bf Y} : \mathbf{Y} \to \mathbb{R}$ is measurable along the flow and the integral
$$
\int_{(0,+\infty)} e^{-\alpha t} C^g_{\bf Y}(\phi_{\bf Y}(y,t))\, dt
$$
is finite.

The key  notations of ${\cal F}^V_+(x)$ and $\underline{\cal F}^V_-(x)$ transform to
$$
\mathcal{F}_+^V(y) \defi \lim_{t \to 0^+} \Big[ \frac{e^{-\alpha t} V(\phi_{\bf Y}(y,t)) - V(y)}{t} + \frac{1}{t} \int_{(0,t]} e^{-\alpha u} C^g_{\bf Y}(\phi_{\bf Y}(y,u))\, du \Big];
$$
$$
\underline{\mathcal{F}}_-^V(y) \defi  \bigg\{ \underline{\lim}_{t \to 0^+} \bigg[ \frac{V(y) - e^{\alpha t} V(\phi_{\bf Y}(\tilde{y}, s-t))}{t} + \frac{1}{t} \int_{(-t,0]}e^{-\alpha u} C^g_{\bf Y}(\phi_{\bf Y}(\tilde{y}, s+u))\, du \Big]:$$
$$(\tilde y,s)\subset {\bf Y}\times \RR_+,~\phi_{\bf Y}(\tilde y,s)=y \bigg\}.
$$

If $\mathbf{Y}'$ is a smooth open manifold, the flow is given by the differential equation $\dot y = f(y)$, satisfying the standard conditions on the existence of a unique local solution in ${\bf Y}'$ (for positive and negative $t$) for each initial condition from ${\bf Y}'$,
{and $C^g_{\bf Y}(y)$ is continuous along the flow in ${\bf Y}'$} and $V$ is continuously differentiable along the flow in ${\bf Y}'$ then $\underline{\mathcal{F}}^V_-(y)$ is a singleton for all $y\in{\bf Y}'$ and
$$
\mathcal{F}_+^V(y) = \underline{\mathcal{F}}_-^V(y) = -\alpha V(y) + C^g_{\bf Y}(y) + \nabla V(y) \cdot f(y).
$$

The integral Bellman equation (\ref{IntEq}) takes the form
$$
V(y) = \inf_{\theta \in \bar{\mathbb{R}}_+^0} \Big\{ \int_{(0,\theta]} e^{-\alpha u} C^g_{\bf Y}(\phi_{\bf Y}(y,u))\, du
+ I\{ \theta < +\infty \}\, e^{-\alpha\theta} \mathcal{I}V(\phi_{\bf Y}(y, \theta)) \Big\},
$$
$$
\text{where $\mathcal{I}V$ is as before, \hspace{16mm}} \mathcal{I}V(y) =  \inf_{a \in \mathbf{A}} \{ C^I_{\bf Y}(y, a) + V(l_{\bf Y}(y, a)) \}.
$$
To be more precise, one had to denote the above Bellman function as $V_{\bf Y}$. In the framework of the extended state space ${\bf X}={\bf Y}\times\RR^0_+$, the Bellman function {is $V((y,s))=e^{-\alpha s}V_{\bf Y}(y)$.}

The Bellman equation in the differential form (\ref{DifEq1}) remains as it was with the obvious changes $x\to y$, $C^I(x,a)\to C^I_{\bf Y}(y,a)$ etc.

Finally, Theorem \ref{th1} remains {valid}, provided that the integral $\int_{(0,+\infty)} e^{-\alpha t} V(\phi_{\bf Y}(y,t))\, dt$ is finite.
{The latter holds true if $C^g_{\bf Y}(y)\le K<\infty$: Condition \ref{con2} is satisfied because $C^g((y,s))= e^{-\alpha s} C^g_{\bf Y}(y)$ and, according to Theorem \ref{t1}, the minimal positive solution to equation (\ref{PZZeqn01}) (i.e., to equation (\ref{IntEq})) has the form $V((y,s))=e^{-\alpha s}V_{\bf Y}(y)$, where $0\le V_{\bf Y}(y)\le \frac{K}{\alpha}$.

As was mentioned after Remark \ref{rem5}, our model is the special case of PDMP when $\lambda=0$. In this connection, it is worth comparing equation $\mathcal{F}_+^V(y) =0$, which comes to the stage only if the right limit exists, and the corresponding differential forms obtained in \cite{b3,b5,b7,b12,b9}. After adding and subtracting $e^{-\alpha t} V(y)$ in the formula for $\mathcal{F}_+^V(y)$, we obtain
$$
\mathcal{F}_+^V(y) = -\alpha V(y) + \lim_{t \to 0^+} \Big[ \frac{e^{-\alpha t} V(\phi_{\bf Y}(y,t)) - e^{-\alpha t} V(y)+ \int_{(0,t]} e^{-\alpha u} C^g_{\bf Y}(\phi_{\bf Y}(y,u))\, du}{t} \Big].
$$
After denoting
$$
\mathcal{X}V(y) \defi  \lim_{t \to 0^+} \Big[ \frac{e^{-\alpha t} V(\phi_{\bf Y}(y,t)) - e^{-\alpha t} V(y)+ \int_{(0,t]} e^{-\alpha u} C^g_{\bf Y}(\phi_{\bf Y}(y,u))\, du}{t} \Big]-C^g_{\bf Y}(y),
$$
equality $\mathcal{F}_+^V(y) =0$ takes the form
$$\mathcal{X}V(y)-\alpha V(y)+C^g_{\bf Y}(y)=0,$$
which appears in \cite{b3,b5,b7,b12,b9} for the case ${\bf X}\subseteq\RR^d$. Moreover, in the smooth case, if the flow comes from the differential equation $\dot y=f(y)$, as was mentioned in Remark \ref{rem2}, ${\cal X} V(y)=\nabla V(y) \cdot f(y)$. (See \cite{b5,b12,b9}.)

Connection between the integral and differential forms of the optimality equation was underlined in \cite{b3,b7,b9}. But it seems that the formal rigorous equivalence of such representations, as established in Theorem \ref{th1}, is proved here for the first time for a general Borel space $\bf X$ and both for discounted and undiscounted cases. As explained in Conclusion, one can easily generalize Theorem \ref{th1} for PDMP.}

\section{{Impulse} Control of SIR Epidemic}
\label{sec:SIR}
In this section, we illustrate all the theoretical issues on a meaningful example having its own interest.

In the following proposition, we enlist all the conditions, which appeared in the previous sections, needed for the study of the model stated below.

\begin{proposition}\label{prqm}
Suppose Conditions \ref{asa}, \ref{PPZcondition01} and \ref{as3} are satisfied. Assume that a lower semicontinuous bounded function $V:~{\bf X}\to\RR^0_+$ is such that
\begin{itemize}
\item equation (\ref{DifEq1}) is valid;
\item Conditions \ref{con4}, \ref{con5} and \ref{con6} are satisfied;
\item inequality $\int_{(0,\infty)} V(\phi(x,t))dt<\infty$ holds true for all $x\in{\bf X}$;
\item $\lim_{i\to\infty} E^\pi_{x_0}[V(X_i)]=0$ for all strategies $\pi$ and for all $x_0\in{\bf X}$.
\end{itemize}
Then
\begin{itemize}
\item $V={\cal V}^*$ and
\item the strategy $(\varphi^*_\theta,\varphi^*_a)$, such that $\varphi^*_\theta(x)=\inf\{\theta:~\phi(x,\theta)\in{\cal L}\}$, where the set ${\cal L}$ is defined in (\ref{e29}), and $\varphi^*_a$ satisfies equality $C^I(x,\varphi^*_a(x))+V(l(x,\varphi^*_a(x)))-V(x)=0$, is uniformly optimal, provided that the both maps $\varphi^*_\theta$ and $\varphi^*_a$ are measurable.
\end{itemize}
\end{proposition}

To formulate our Susceptible--Infected--Recovered (SIR) model
of epidemics, we use functions $t \mapsto x^1(t)$,
$t \mapsto x^2(t)$ and $t \mapsto x^3(t)$,
where $x^1 : {\mathbb{R}}^0_+ \rightarrow \mathbb{R}^0_+$
denotes the dynamics of the susceptible population,
$x^2 : {\mathbb{R}}^0_+ \rightarrow \mathbb{R}^0_+$
the dynamics of the infective population  and
$x^3 : {\mathbb{R}}^0_+ \rightarrow \mathbb{R}^0_+$
the dynamics of the removed population (recovered or dead).
Following \cite{MR0452809,MR2199469,MR0945530,MR2478635}, the progress
of infection is described by the following initial value problem:
\begin{equation}
\label{eq:ivp}
\begin{gathered}
\begin{cases}
\dot{x}^1(t) = -\beta \displaystyle \frac{x^1(t) x^2(t)}{x^1(t)+x^2(t)},\\[0.3cm]
\dot{x}^2(t) = \beta \displaystyle \frac{x^1(t) x^2(t)}{x^1(t)+x^2(t)} - \gamma x^2(t),\\[0.3cm]
\dot{x}^3(t) = \gamma x^2(t),
\end{cases}\\
x^1(0) = x^1_0, \quad x^2(0) = x^2_0, \quad x^3(0) = 0,
\end{gathered}
\end{equation}
for some given constant parameters $\beta, \gamma \in \mathbb{R}_+$
and initial data $x^1_0, x^2_0 \in \mathbb{R}^0_+$. If $x^1_0=x^2_0=0$ then $x^1(t)\equiv x^2(t)\equiv 0$.
Problem \eqref{eq:ivp} has a unique solution, obtained in the closed
form in \cite{MR0452809,MR0945530}. Explicit (non-impulse) optimal control policies
for \eqref{eq:ivp} are available in the literature: optimal isolation/treatment
of infective individuals has been studied in \cite{MR2199469} while the case
of immunization/vaccination is investigated in \cite{MR2478635}. Here we
formulate and explicitly solve an optimal control problem with
isolation/treatment impulses.


\subsection{Optimal Control Problem}
\label{sec:SIR:OCP}

Suppose there are no impulses.

We begin by noting that in \eqref{eq:ivp} one has
$\dot{x}^1(t) + \dot{x}^2(t) + \dot{x}^3(t) = 0$ for any $t$,
so that the total population is constant along time:
$x^1(t) + x^2(t) + x^3(t)$ is a fixed constant $x^1_0 + x^2_0 \in \mathbb{R}_+$.
For this reason, $x^3(t) = x^1_0+x^2_0 - x^1(t) - x^2(t)$ and it is sufficient to
restrict ourselves to differential equations
\begin{equation}
\label{eq:ivp:xy}
\begin{gathered}
\begin{cases}
\dot{x}^1(t) = -\beta \displaystyle \frac{x^1(t) x^2(t)}{x^1(t)+x^2(t)},\\[0.3cm]
\dot{x}^2(t) = \beta \displaystyle \frac{x^1(t) x^2(t)}{x^1(t)+x^2(t)} - \gamma x^2(t),
\end{cases}\\
x^1(0) = x^1_0, \quad x^2(0) = x^2_0,
\end{gathered}
\end{equation}
which define the flow $\phi$.

Since there is no immigration (and no births) and isolation leads to the decrease of $x^2$, the whole state space is the triangle
$${\bf X}=\{x^1\ge 0,~x^2\ge 0,~~x^1+x^2 < N \}$$
with the topology induced from $\RR^2$;
it is convenient to take $N>x^1_0+x^2_0$, so that there are no singular points in $\bf X$.
The gradual cost rate is the infection rate
\begin{equation}
\label{eq:SIR:GCR}
C^g(x^1,x^2) = \beta \frac{x^1 x^2}{x^1+x^2},~~~C^g(0,0)=0,
\end{equation}
which, after integration along the flow, results in the total number of new infectives. Here and below, usually the brackets of the argument of a function of $x=(x^1,x^2)$ are omitted.

At any moment, the decision maker can isolate all infectives, so that ${\bf A}=\{1\}$ is a singleton. The cost of an impulse equals
\begin{equation}
\label{eq:SIR:CI}
C^I(x^1,x^2,a) = c x^2,
\end{equation}
where $c > 0$ is a given constant. The new state after the impulse equals
\begin{equation}
\label{eq:SIR:NSAI}
l(x^1,x^2,a) = (x^1, 0).
\end{equation}
If $x^2_0=0$ then
 there are no individuals who can cause infection
and, therefore, the susceptible population will remain constant forever:
$x^1(t) \equiv x^1_0$. Thus
$${\bf Y}=\{(x^1,x^2)\in{\bf X}:~x^2=0\}$$
is the "cemetery" subset.

\begin{remark}\label{rem82}
Quite formally, in the states $(x^1,x^2)\in{\bf Y}$, one can still apply impulses leading to no cost and no change of the state. But actually, the controlled process is finished as soon as the state belongs to $\bf Y$.
\end{remark}

Note also that
\begin{equation}\label{e39}
\int_{(0,+\infty)} C^g(\phi((x^1_0,x^2_0),u))du \le x^1_0< { N <} +\infty,
\end{equation}
meaning that the Bellman function ${\cal V}^*$ is bounded on the bounded domain $\bf X$.

One can easily check that  Conditions \ref{asa},\ref{PPZcondition01},\ref{con2} and \ref{as3} are satisfied.

To solve the optimal control problem, we investigate the differential form of the Bellman equation (\ref{DifEq1}) which is equivalent to (\ref{PZZeqn01}) by Theorem \ref{th1}. As was mentioned in Remark \ref{rem2}, under certain conditions { which are satisfied in our example},
\begin{equation}\label{e23}
{\cal F}^V_+(x^1,x^2) = \underline{\cal F}^V_-(x^1,x^2)= { \beta\frac{x^1x^2}{x^1+x^2}} +\frac{\partial V}{\partial x^1}\left(-\beta\frac{x^1x^2}{x^1+x^2}\right)+\frac{\partial V}{\partial x^2}\left(\beta\frac{x^1x^2}{x^1+x^2}-\gamma x^2\right).
\end{equation}

In the future, we will need the following explicit expressions defining the flow $\phi(x^1,x^2)$ at $x^1,x^2 > 0$ coming from the differential equation (\ref{eq:ivp:xy}) (see \cite{MR0945530,MR2478635}):
\begin{equation}\label{e24}
\left.\begin{array}{ll}
\mbox{if } \beta\ne\gamma, & \mbox{ then} \\
& \displaystyle
x^1(t)=x^1_0\frac{\left(1+\frac{x^2_0}{x^1_0}\right)^{\frac{\beta}{\beta-\gamma}}}{\left(1+\frac{x^2_0}{x^1_0} e^{(\beta-\gamma)t}\right)^{\frac{\beta}{\beta-\gamma}}};\\
& \displaystyle
x^2(t)=x^2_0\frac{\left(1+\frac{x^2_0}{x^1_0}\right)^{\frac{\beta}{\beta-\gamma}} e^{(\beta-\gamma)t}}{\left(1+\frac{x^2_0}{x^1_0} e^{(\beta-\gamma)t}\right)^{\frac{\beta}{\beta-\gamma}}};\\
\mbox{if } \beta=\gamma, & \mbox{ then} \\
& \displaystyle
x^1(t)=x^1_0 e^{-\frac{\beta x^2_0 t}{x^1_0+x^2_0}};~~~x^2(t)=x^2_0 e^{-\frac{\beta x^2_0 t}{x^1_0+x^2_0}}.
\end{array}\right\}
\end{equation}
From these expressions, it is clear that
\begin{equation}\label{e31}
\frac{x^2(t)}{x^1(t)}=\frac{x^2_0}{x^1_0} e^{(\beta-\gamma)t}
\end{equation}
and, for all $t\ge 0$, $x^1(t),x^2(t)>0$, if $x^1_0,x^2_0>0$.

{ Below, we summarise general properties of the epidemic model under study.
\begin{itemize}
\item The Bellman function ${\cal V}^*$ is bounded and Condition \ref{as3} is satisfied because of inequality (\ref{e39}).
\item Conditions \ref{asa} and \ref{PPZcondition01} are satisfied.
\item If $V$ is a bounded lower semicontinuous function satisfying equation (\ref{PZZeqn01}), then $V(0,x^2)=V(x^1,0)=0$ for all $x^1,x^2\in\RR^0_+$. As a result, for all strategies $\pi$, for all $(x^1_0,x^2_0)\in{\bf X}$, $E^\pi_{(x^1_0,x^2_0)}[V(X_1)]=0$ because either $X^2_1=0$ (if $\Theta_1<\infty$) and all the further values of the second component equal zero, or the next state is $X_1=\Delta$ (if $\Theta_1=\infty$). Hence, $\lim_{i\to\infty} E^\pi_{x_0}[V(X_i)]=0$ for all strategies $\pi$ and for all $x_0\in{\bf X}$.
\item According to Remark \ref{rem9}, the Bellman equation (\ref{PZZeqn01}) has a unique bounded lower semicontinuous solution $V={\cal V}^*$.
\item If function $V:~{\bf X}\to\RR$ is such that $V(x^1,0)=V(0,x^2)=0$ then, for $x^1=0$ or $x^2=0$, equalities
$${\cal F}^V_+(x^1,x^2)=\underline{\cal F}^V_-(x^1,x^2)=0$$
and
$$\inf_{a\in{\bf A}}[C^I(x^1,x^2,a)+V(l(x^1,x^2,a))-V(x^1,x^2)]=cx^2$$
are valid. Thus, equations (\ref{DifEq1}) hold: case (a) is $x^2>0$ and $x^1=0$ and case (b) if $x^2=0$. All the conditions \ref{con4}, \ref{con5} and \ref{con6} are trivially satisfied if $x^1=0$ or $x^2=0$.
\end{itemize}

The points of the form  $(x^1,0)$ belong to $\cal L$. Formally speaking, if $x^2=0$ one has to apply the simultaneous infinite sequence of impulses, each of them having no effect: see Remark \ref{rem82}. At the states $(0,x^2>0)$, no impulses are needed: the number of infectives $x^2$ decreases to zero resulting in no cost.

We have seen that the flow $\phi(x,t)$ is continuous. We will see in all the three cases investigated in the further subsections that the function $V$ to be defined below is continuous on ${\bf X}\cap(\RR_+)^2$ and the corresponding set $\cal L$ is closed. It follows that conditions \ref{con4}, \ref{con5} and \ref{con6} are satisfied.
}
\subsection{Solution in the Case $\beta\ge\gamma$}
\label{sec:SIR:solOCP}

In this subsection, we show that the { continuous} function
$$  V(x^1,x^2)=\left\{\begin{array}{ll}
cx^2, & \mbox{ if } x^2\le\frac{x^1}{c}, ~ x^1 \ge 0, ~ x^2 \ge 0; \\
x^1, & \mbox{ if } x^2>\frac{x^1}{c}, ~ x^1 \ge 0,~x^2>0
\end{array}\right. $$
satisfies all the requirements of { Proposition \ref{prqm}}.

Firstly, let us show that the integral $\int_{(0,+\infty)} V(\phi(x^1_0,x^2_0,t))dt$ is finite for all $x^1_0,x^2_0>0$.

If $x^2_0 > \frac{x^1_0}{c}$, then
$$x^2(t) > \frac{x^1(t)}{c} \mbox{ for all }t > 0$$
because of (\ref{e31}).
Therefore, for such initial values $(x^1_0,x^2_0)$,
$$\int_{(0,+\infty)} V(\phi(x^1_0,x^2_0,t))dt=\int_{(0,+\infty)} x^1(t)dt.$$
According to (\ref{e31}),
$$x^1(t)=x^2(t)\frac{x^1_0}{x^2_0} e^{-(\beta-\gamma)t}.$$
Since function $x^2(t)$ is uniformly bounded, in case $\beta>\gamma$, the integral $\int_{(0,+\infty)} V(\phi(x^1_0,x^2_0,t))dt$ is finite. If $\beta=\gamma$, its finiteness follows directly from (\ref{e24}).

In case $x^2_0 \le \frac{x^1_0}{c}$ and $\beta>\gamma$, again using (\ref{e31}), we see that $x^2(t^*)=\frac{x^1(t^*)}{c}$ at $t^*=\frac{\ln(x^1_0)-\ln(x^2_0 c)}{\beta-\gamma}<\infty$ and
$$\int_{(0,+\infty)} V(\phi(x^1_0,x^2_0,t))dt=\int_{(0,t^*]} c x^2(t) dt+\int_{(t^*,+\infty)} x^1(t)dt<\infty.$$
If $x^2_0 \le \frac{x^1_0}{c}$ and $\beta=\gamma$, then $x^2(t)\le \frac{x^1(t)}{c}$ for all $t\ge 0$ and $\displaystyle\int_{(0,+\infty)} V(\phi(x^1_0,x^2_0,t))dt<\infty$ by (\ref{e24}).

The { closed} set $\cal L$ defined in (\ref{e29}) has the form
$${\cal L}=\big\{(x^1,x^2)\in{\bf X}\cap(\RR^0_+)^2:~x^2\le\frac{x^1}{c}\big\}.$$

Now show that equation (\ref{DifEq1}) is valid.

If $(x^1,x^2)\in{\cal L}$ then
\begin{eqnarray*}
\underline{\cal F}^V_-(x^1,x^2)&=&\beta\frac{x^1 x^2}{x^1+x^2}+\frac{\partial V}{\partial x^2}\left[\beta\frac{x^1 x^2}{x^1+x^2}-\gamma x^2\right]
= \frac{x^2}{x^1+x^2}\left[ \beta(1+c)x^1-\gamma c(x^1+x^2)\right].
\end{eqnarray*}
On the boundary $x^2=\frac{x^1}{c}$, the expression $\frac{\partial V}{\partial x^2}$ means the left derivative, and the {right} derivative $\frac{\partial V}{\partial x^1}=0$.
Since $x^2\le\frac{x^1}{c}$, we conclude that
$$\underline{\cal F}^V_-(x^1,x^2)\ge \frac{x^1 x^2}{x^1+x^2}[\beta(1+c)-\gamma(c+1)]=\frac{x^1 x^2}{x^1+x^2} (\beta-\gamma)(c+1)\ge 0,$$
so that equality (\ref{DifEq1}), case (b), is satisfied.

The cases $x^1=0$ or $x^2=0$ were considered in Subsection \ref{sec:SIR:OCP}.

If $(x^1,x^2)\notin{\cal L}$ and $(x^1,x^2)\in(\RR_+)^2$ then
$$\inf_{a\in{\bf A}} [C^I(x^1,x^2,a)+V(l(x^1,x^2,a))-V(x^1,x^2)]>0$$
and
$${\cal F}^V_+(x^1,x^2)=\beta\frac{x^1 x^2}{x^1+x^2} +\frac{\partial V}{\partial x^1}\left[-\beta\frac{x^1 x^2}{x^1+x^2}\right]=0,$$
so that equality (\ref{DifEq1}), case (a), is satisfied.

According to { Proposition \ref{prqm}}, the stationary strategy
$$\varphi^*_\theta(x^1,x^2) = \left\{\begin{array}{lll}
\infty, & \mbox{ if } (x^1,x^2)\notin{\cal L} & \Longleftrightarrow x^2>\frac{x^1}{c};\\
0, & \mbox{ if } (x^1,x^2)\in{\cal L} & \Longleftrightarrow  x^2\le \frac{x^1}{c};
\end{array}\right.~~~\varphi^*_a(x^1,x^2)=1 $$
is uniformly optimal. {The straight line $x^2=\frac{x^1}{c}$ is a dispersal line.}

For $(x^1,x^2)\in(\RR_+)^2$, it is reasonable to rewrite expression $(x^1,x^2)\in{\cal L}$ as
$$\{(x^1,x^2)\in(\RR_+)^2:~~x^1\ge c x^2\},$$
to understand better the meaning of the optimal strategy. The goal of the control is to save susceptibles  from being infected, but the cost of isolation is $c x^2$. Thus, isolation is reasonable only when there are many susceptibles to be saved: $x^1\ge c x^2$ because otherwise the cost of isolation, $c x^2$, is bigger than the profit for saving susceptibles (i.e., $x^1$).

\begin{figure}[!htb]
	\centering
	\includegraphics[scale=0.7]{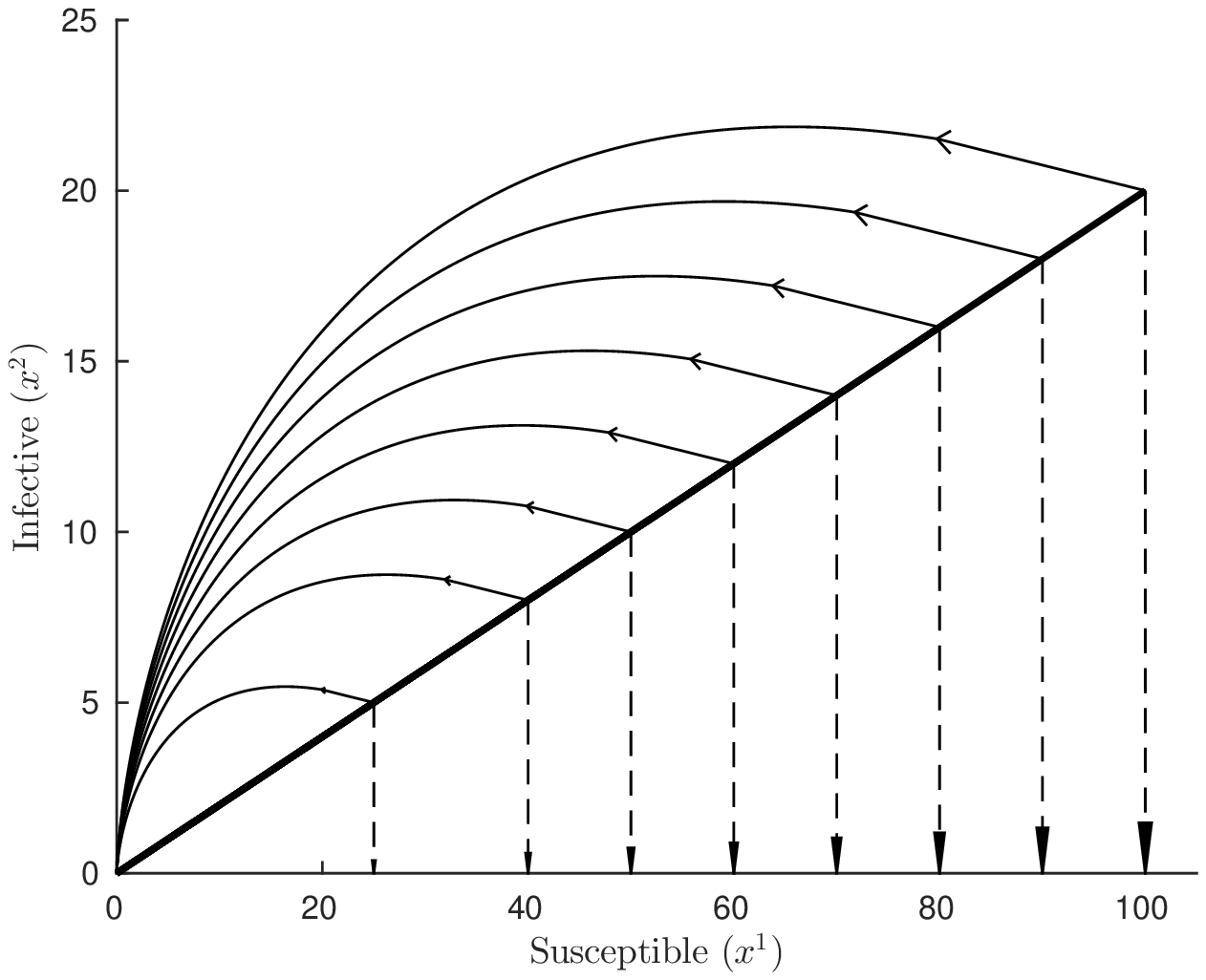}
	\caption{Susceptible--Infected dynamics under optimal control
		with $c = 5$, $\beta = 4$ and $\gamma = 3$.}
	\label{fig:01}

	\centering
	\includegraphics[scale=0.7]{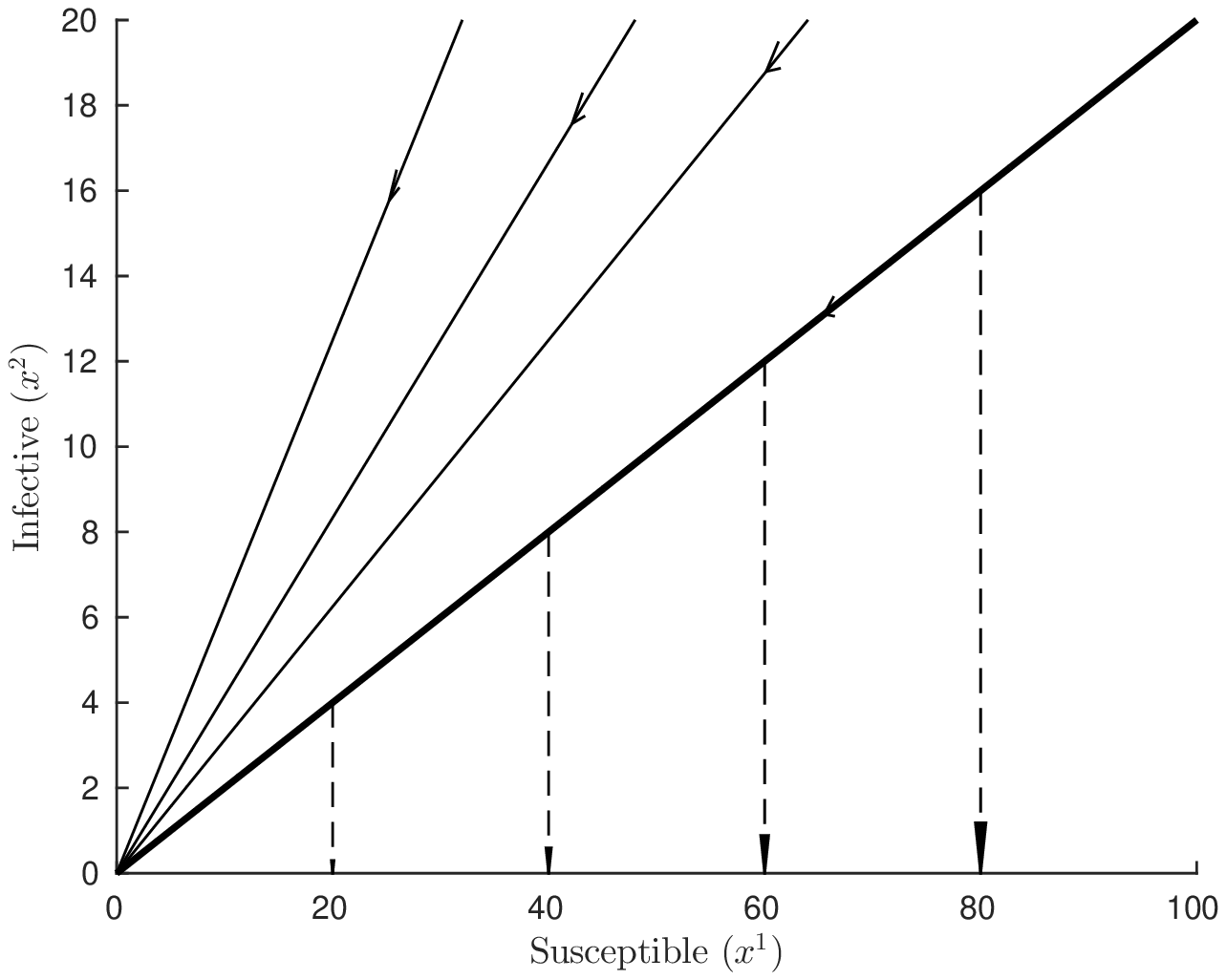}
	\caption{Susceptible--Infected dynamics under optimal control
		with $c = 5$, $\beta = 4$ and $\gamma = 4$.}
	\label{fig:04}
\end{figure}

The optimal strategy is shown in Figures \ref{fig:01} and \ref{fig:04}. If the initial state $(x_0^1,x_0^2)$ is below the line $x^2=\frac{x^1}{c}$ (shown in bold) then the impulse should be applied (dashed line).  If, otherwise, the initial state is above the line $x^2=\frac{x^1}{c}$ then no impulse is needed and the system evolves according to equations \eqref{e24} (solid curves). If $\beta=\gamma$, then the critical line $x^2=\frac{x^1}{c}$ is the trajectory of the dynamical system (\ref{e24}). It is equally optimal to move along this line or to apply the impulse immediately or at any further time.

\subsection{Solution in the Case $\beta<\gamma$}
\label{subsec3}
\subsubsection{Case $c\ge\frac{\beta}{\gamma-\beta}$} \label{s531} In this subsection, we show that the function
$$V(x^1,x^2)=x^1-x^1\left(1+\frac{x^2}{x^1}\right)^{-\frac{\beta}{\gamma-\beta}},~~~(x^1,x^2)\in(\RR_+)^2$$
satisfies all the requirements of { Proposition \ref{prqm}}. According to Subsection \ref{sec:SIR:OCP}, $V(x^1,x^2)=0$ if $x^1=0$ or $x^2=0$.
Firstly, let us show that the integral $\int_{(0,+\infty)} V(\phi(x^1_0,x^2_0,t))dt$ is finite for all $x^1_0,x^2_0>0$. According to (\ref{e31}) and keeping in mind that $x^1(t)\le x^1_0$, it is sufficient to prove that the integral
$$I=\int_{(0,+\infty)} \left[1-\left(1+\frac{x^2_0}{x^1_0} e^{-(\gamma-\beta)t}\right)^{-\frac{\beta}{\gamma-\beta}}\right] dt$$
is finite. { This is a simple consequence of the fact that the integrand  is $O(e^{-(\gamma-\beta)t})$ as $t\to\infty$.}

In the case under consideration, ${\cal L}={\bf Y}=\{(x^1,x^2):~x^2=0\}$ { is closed}. The value $x^2=0$ is not reachable in finite time from initial conditions $(x^1,x^2>0)$.

It remains to check equation (\ref{DifEq1}) for the presented function $V$. Namely, we will show that the version (a) is valid. The cases $x^1=0$ or $x^2=0$ were considered in Subsection \ref{sec:SIR:OCP}. For $(x^1,x^2)\in(\RR_+)^2$, according to (\ref{e23}), equality ${\cal F}^V_+(x^1,x^2)=0$ can be checked straightforwardly. Finally, for $(x^1,x^2)\in(\RR_+)^2$,
$$\inf_{a\in{\bf A}} [C^I(x^1,x^2,a)+V(l(x^1,x^2,a))-V(x^1,x^2)]=cx^2-x^1+x^1\left(1+\frac{x^2}{x^1}\right)^{-\frac{\beta}{\gamma-\beta}} $$
$$
\ge x^1 \left[ \left(1+\frac{x^2}{x^1}\right)^{-\frac{\beta}{\gamma-\beta}} - \left( 1 - \frac{\beta}{\gamma-\beta}\frac{x^2}{x^1} \right)\right] >0.$$

According to { Proposition \ref{prqm}}, the stationary strategy
$$\varphi^*_\theta(x^1,x^2) = \left\{\begin{array}{lll}
\infty, & \mbox{ if } (x^1,x^2)\notin{\cal L} & \Longleftrightarrow x^2>0;\\
0, & \mbox{ if } (x^1,x^2)\in{\cal L} & \Longleftrightarrow  x^2=0;
\end{array}\right.~~~\varphi^*_a(x^1,x^2)=1 $$
is uniformly optimal.

When $x^2>0$, isolation is not reasonable as its cost $c\ge \frac{\gamma}{\gamma-\beta}$ is too high. When $x^2=0$, the epidemic is actually terminated, although the formal solution prescribes isolation of zero infectives for zero cost, without any real effect.

\begin{figure}[!htb]
	\centering
	\includegraphics[scale=0.7]{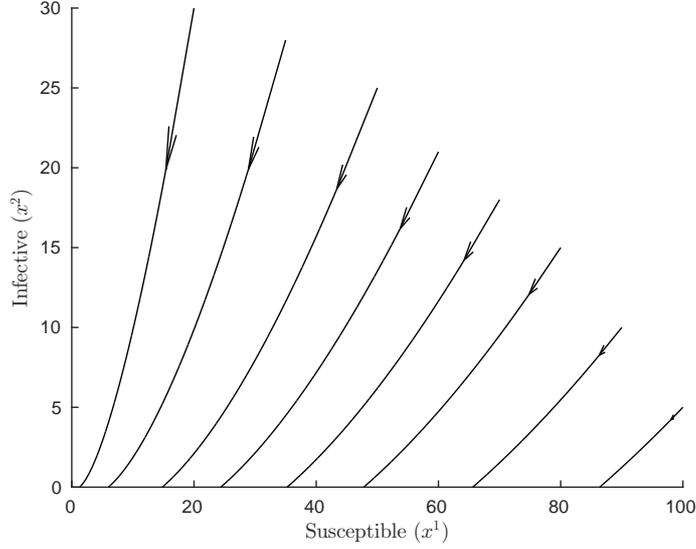}
	\caption{Susceptible--Infected dynamics with $c = 5$, $\beta = 3$ and $\gamma = 4$.}
	\label{fig:02}
\end{figure}

The optimal strategy in this case for the values $c = 5$, $\beta = 3$ and $\gamma = 4$ is shown in Figure~\ref{fig:02}. No impulses are needed here, and the system evolves according to equations \eqref{e24}.


\subsubsection{Case $c<\frac{\beta}{\gamma-\beta}$}
\label{sebsec4}

In this subsection, we show that the continuous function
$$
V(x^1,x^2)=\left\{\begin{array}{ll}
c x^2, & \mbox{ if } {x^2}\le \frac{\beta+\beta c-\gamma c}{\gamma c}\,{x^1};\\
x^1\left[ 1-\left(\frac{\gamma c\left(1+\frac{x^2}{x^1}\right)}{\beta+\beta c}\right)^{-\frac{\beta}{\gamma-\beta}}
(1+c)\,\frac{\gamma-\beta}{\gamma}
\right], & \mbox{ if } {x^2}>\frac{\beta+\beta c-\gamma c}{\gamma c}\,{x^1}
\end{array}\right.$$
satisfies all the requirements of { Proposition \ref{prqm}}.

Firstly, let us show that the integral $\int_{(0,+\infty)} V(\phi(x^1_0,x^2_0,t))dt$ is finite for all $x^1_0,x^2_0>0$. Indeed, if $\frac{x^2_0}{x^1_0}\le\frac{\beta+\beta c-\gamma c}{\gamma c}$, then, by (\ref{e31}),
$$\int_{(0,+\infty)} V(\phi(x^1_0,x^2_0,t))dt=c\int_{(0,+\infty)} x^1(t)\left(\frac{x^2_0}{x^1_0} e^{-(\gamma-\beta)t}\right) dt<\infty$$
because the function $x^1(t)\le x^1_0$ is bounded. If $\frac{x^2_0}{x^1_0}>\frac{\beta+\beta c-\gamma c}{\gamma c}$, then $\frac{x^2(t^*)}{x^1(t^*)}=\frac{\beta+\beta c-\gamma c}{\gamma c}$ at the finite time moment
$$t^*=\theta^*(x^1,x^2)=\frac{1}{\gamma-\beta} \ln \frac{\gamma c x^2_0}{x^1_0(\beta+\beta c-\gamma c)}>0,$$
the integral $\int_{(0,t^*]} V(\phi(x^1_0,x^2_0,t))dt$ is finite, and on the interval $(t^*,+\infty)$ the previous reasoning applies.

{ Let us check that the set $\cal L$ defined in (\ref{e29}) has the form}
$${\cal L}=\{(x^1,x^2)\in{\bf X}\cap(\RR^0_+)^2:~x^2\le \frac{\beta+\beta c-\gamma c}{\gamma c}\cdot x^1\},$$
{ and therefore is closed.}

The case when $x^1=0$ or $x^2=0$ was considered in Subsection \ref{sec:SIR:OCP}.

If $0<x^2\le \frac{\beta+\beta c-\gamma c}{\gamma c}\cdot x^1$, then
$$C^I(x^1,x^2,a)+V(l(x^1,x^2,a))-V(x^1,x^2)= cx^2+V(x^1,0)-cx^2=0,$$
so that $(x^1,x^2)\in{\cal L}$. Remember, $a=1\in{\bf A}$ is the unique action.

If $x^2 > \frac{\beta+\beta c-\gamma c}{\gamma c}\cdot x^1>0$, then
$$C^I(x^1,x^2,a)+V(l(x^1,x^2,a))-V(x^1,x^2)= x^1\Upsilon\left(\frac{x^2}{x^1}\right),$$
where { function}
$$\Upsilon(w)=cw-1+\left(\frac{\gamma c\left(1+w\right)}{\beta+\beta c}\right)^{-\frac{\beta}{\gamma-\beta}}
 (1+c)\,\frac{\gamma-\beta}{\beta}$$
{ is strictly convex.} When $w=\frac{\beta+\beta c-\gamma c}{\gamma c}$,
$$\Upsilon\left(\frac{\beta+\beta c-\gamma c}{\gamma c}\right)=\frac{d\Upsilon}{dw}\left(\frac{\beta+\beta c-\gamma c}{\gamma c}\right)=0.$$
Thus, $\Upsilon (w) > 0$ for  $w>\frac{\beta+\beta c-\gamma c}{\gamma c}$, and
therefore,
$$C^I(x^1,x^2,a)+V(l(x^1,x^2,a))-V(x^1,x^2)>0$$
if $x^2>\frac{\beta+\beta c-\gamma c}{\gamma c}\cdot x^1\ge 0$, and $(x^1,x^2)\notin{\cal L}$: the case $x^1=0$ is not excluded, as well.

Now show that equation (\ref{DifEq1}) is valid.

If $(x^1,x^2)\notin{\cal L}$ and $(x^1,x^2)\in(\RR_+)^2$ then we already know that
$$C^I(x^1,x^2,a)+V(l(x^1,x^2,a))-V(x^1,x^2)>0.$$
Equality
\begin{equation}\label{e333}
{\cal F}^V_+(x^1,x^2)=\beta\frac{x^1 x^2}{x^1+x^2} +\frac{\partial V}{\partial x^1}\left[-\beta\frac{x^1 x^2}{x^1+x^2}\right]
+\frac{\partial V}{\partial x^2}\left[\beta\frac{x^1 x^2}{x^1+x^2}-\gamma x^2\right]=0
\end{equation}
in the area $x^2>\frac{\beta+\beta c-\gamma c}{\gamma c}\cdot x^1>0$ can be checked by the direct substitution.
Equation (\ref{DifEq1}), case (a), is valid.

The cases $x^1=0$ or $x^2=0$ were considered in Subsection \ref{sec:SIR:OCP}.

If $x^1>\frac{\gamma c}{\beta+\beta c-\gamma c}\,x^2>0$, then
\begin{eqnarray*}
\underline{\cal F}^V_-(x^1,x^2)&=&\beta\frac{x^1 x^2}{x^1+x^2}+\frac{\partial V}{\partial x^2}\left[\beta\frac{x^1 x^2}{x^1+x^2}-\gamma x^2\right]
= \frac{x^2}{x^1+x^2}\left[ \beta(1+c)x^1-\gamma c(x^1+x^2)\right]\\
&\ge & \frac{x^2}{x^1+x^2}\left[(\beta+\beta c-\gamma c)\frac{x^2\gamma c}{\beta+\beta c-\gamma c}-\gamma c x^2 \right]=0.
\end{eqnarray*}
On the boundary $x^2=\frac{\beta+\beta c-\gamma c}{\gamma c}\cdot x^1$, we need to consider the  left derivative $\frac{\partial V}{\partial x^1}$ and the right derivative $\frac{\partial V}{\partial x^2}$.  As  a result, $\underline{\cal F}^V_-(x^1,x^2)=0$ similarly to (\ref{e333}).
Equation (\ref{DifEq1}), case (b), is valid.

According to Proposition \ref{prqm}, the stationary strategy
$$\varphi^*_\theta(x^1,x^2) = \left\{\begin{array}{lll}
\infty, & \mbox{ if } (x^1,x^2)\notin{\cal L} & \Longleftrightarrow x^2>\frac{\beta+\beta c-\gamma c}{\gamma c}\cdot x^1;\\
0, & \mbox{ if } (x^1,x^2)\in{\cal L} & \Longleftrightarrow  x^2\le\frac{\beta+\beta c-\gamma c}{\gamma c}\cdot x^1;
\end{array}\right.~~~\varphi^*_a(x^1,x^2)=1 $$
is uniformly optimal. {The straight line $x^2=\frac{\beta+\beta c-\gamma c}{\gamma c}\cdot x^1$ is a switching line.}

Like in the case $\beta\ge \gamma$, isolation of infectives is reasonable only when there are sufficiently many susceptibles to be saved: $x^1\ge\frac{ \gamma c}{\beta+\beta c-\gamma c}\, x^2$.

\begin{figure}[!htb]
	\centering
	\includegraphics[scale=0.7]{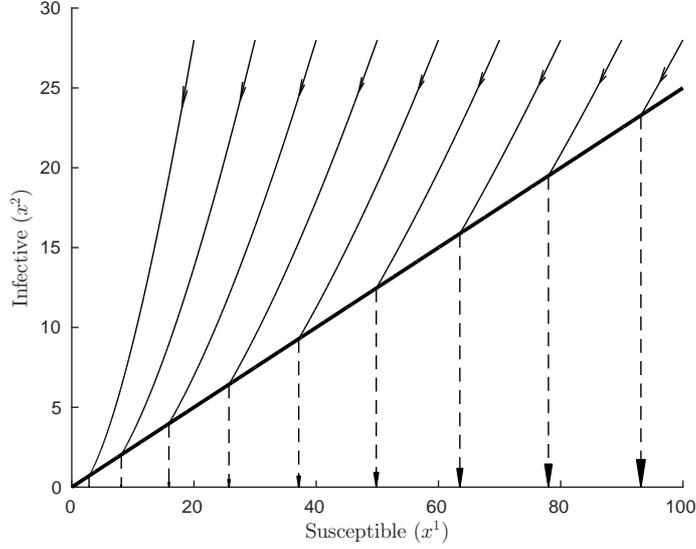}
	\caption{Susceptible--Infected dynamics under optimal control
		with $c = 3/2$, $\beta = 3$ and $\gamma = 4$.}
	\label{fig:03}
\end{figure}

In this case we take $c = 3/2$, $\beta = 3$ and $\gamma = 4$; see Figure~\ref{fig:03}. If the initial state $(x_0^1,x_0^2)$ lies below the line $x^2=\frac{\beta+\beta c-\gamma c}{\gamma c}\cdot x^1$ (shown in bold) then the impulse should be applied (dashed line).  If the initial state lies above this line then {initially} no action is needed and the system evolves according to equations \eqref{e24} (solid curves)
{up to the moment when $x^2(t)=\frac{\beta+\beta c-\gamma c}{\gamma c}\cdot x^1(t)$ when the impulse should be applied.}

{\subsection{Discussion}\label{subsec4}
The threshold nature of the optimal isolation strategy for other epidemic models with similar cost functions was established in \cite{b14,b13}: intervene only if the current number of infectives is below a certain value. Moreover, it was shown that the intervention must be global, i.e., it is better to isolate all infectives at once.

It is interesting to compare the impulse control problem from Subsection \ref{sec:SIR:OCP} with its gradual control analogue investigated in \cite{MR2199469}. Instead of impulses, dynamic control $u(t)\in[0,U]$ appears in the second equation of (\ref{eq:ivp:xy}):
$$\dot{x}^2(t) = \beta \displaystyle \frac{x^1(t) x^2(t)}{x^1(t)+x^2(t)} - \gamma x^2(t)-u(t)x^2(t).$$
Objective functional in \cite{MR2199469}
$${\cal V}(x_0^1,x_0^2,u)=\int_{(0,\infty)} \left(\beta \frac{x^1(t) x^2(t)}{x^1(t)+x^2(t)}+cu(t) x^2(t)\right)dt\to\inf_u$$
has the same meaning as in the current paper: combination of the total number of the new infectives and the total cost of isolation with the weight coefficient $c>0$. Intuitively, the impulse isolation at time moment $t$ means that $u(t)\to\infty$. Thus, look at the optimal strategy obtained in \cite{MR2199469} when $U\to\infty$.
\begin{itemize}
\item If $\beta\ge\gamma$ then one has to apply the maximal rate of isolation $U$ as soon as  $x^2\le \zeta(U) x^1$, where
$$\zeta(U)=\left(\frac{\gamma+U+cU}{cU}\right)^{\frac{\gamma+U-\beta}{\gamma+U}}-1,$$
and the straight line $x^2=\zeta(U) x^1$ is a dispersal line. When $U\to\infty$,
$$\lim_{U\to\infty} \zeta(U)=\frac{1+c}{c}-1=\frac{1}{c},$$
and we finish with exactly the optimal impulse strategy presented in Subsection \ref{sec:SIR:solOCP}.
\item If $\beta<\gamma$ and $c\ge \frac{\beta}{\gamma-\beta}$ then, both in \cite{MR2199469} and in Subsubsection  \ref{s531}, it is optimal not to immunize at all.
\item If $\beta<\gamma$ and $c<\frac{\beta}{\gamma-\beta}$ then one has to apply the maximal rate of isolation $U$ as soon as $x^2\le \xi(U) x^1$, where
    $$\xi(U)=\left(\frac{\beta(\gamma+U+cU)}{c\gamma(\gamma+U-\beta)}\right)^{\frac{\gamma+U-\beta}{\gamma+U}}-1,$$
and the straight line $x^2=\xi(U) x^1$ is a switching line. When $U\to\infty$,
$$\lim_{U\to\infty} \xi(U)=\frac{\beta(1+c)}{c\gamma}-1=\frac{\beta+\beta c-\gamma c}{\gamma c},$$
and we finish with exactly the optimal impulse strategy presented in Subsection \ref{sebsec4}.
\end{itemize}

There are many other sensible optimal control problems in mathematical epidemiology. For example, one can consider immunization of susceptibles. Such problem for the model (\ref{eq:ivp:xy}) was solved in \cite{MR2478635}, but again in the framework of gradual dynamic control, where the term $-u(t)x^1(t)$ appears in the first equation of (\ref{eq:ivp:xy}). No doubt, the impulse version of immunization can also be tackled using the methods developed in the current paper.}

\section{Conclusion}

Application of the MDP methods to the purely deterministic optimal impulse control problem results in the integral optimality equation. After that, a formal analytical proof shows that the integral and differential forms are equivalent. All the theory is illustrated by a meaningful example on the SIR epidemic.

Note that Theorem \ref{th1} remains also valid in the case when the underlying process is a Piecewise Deterministic Markov Process. To be specific, consider the discounted version of the positive model with the state space ${\bf Y}=\RR^d$, the uncontrolled flow, and the uncontrolled fixed jumps intensity $\lambda$. Under the mild  relevant conditions, the integral equation (\ref{IntEq}) was obtained in \cite{b3,b4,b5,b7,b12}; it has the form
\begin{eqnarray}
V_{\bf Y}(y)&=&\inf_{\theta\in\bar\RR^0_+}\left\{\int_{(0,\theta]} e^{-\alpha u}\left[\vphantom{\int_{\bf Y}} C^g_{\bf Y}(\phi_{\bf Y}(y,u))
+\lambda\left[\int_{\bf Y} V_{\bf Y}(z)Q(dz|\phi_{\bf Y}(y,u))-V_{\bf Y}(\phi_{\bf Y}(y,u))\right]\right]\right.du\nonumber \\
&&+\left. \vphantom{\int_{(0,\theta]}} I\{\theta<+\infty\}e^{-\alpha \theta} \inf_{a\in{\bf A}}\{C^I_{\bf Y}(\phi_{\bf Y}(y,\theta),a)+V_{\bf Y}(l_{\bf Y}(\phi_{\bf Y}(y,\theta),\theta))\}\right\},\label{constar}
\end{eqnarray}
where $Q$ is the stochastic kernel describing the distribution after the spontaneous (natural) jumps with intensity $\lambda>0$. Here, we follow the notations introduced for the discounted model in Section \ref{sec3}, which also appeared in Subsection \ref{ssec42}.

Suppose a measurable along the flow function $V_{\bf Y}:{\bf Y}\to\RR$  is the minimal positive solution to equation (\ref{constar}) and satisfies the corresponding discounted version of Condition \ref{con41}. One can show that, if $C^g_{\bf Y}(y)\le K<\infty$, then the integral $\int_{(0,+\infty)} e^{-\alpha t} V_{\bf Y}(\phi_{\bf Y}(y,t)) dt$ is finite for all $y\in{\bf Y}$. Denote
$$D(y)\defi \lambda\left[\int_{\bf Y} V_{\bf Y}(z) Q(dz|y)-V_{\bf Y}(y)\right].$$
Then, by Theorem \ref{th1}, $V_{\bf Y}(y)$ satisfies the discounted version of the differential equation (\ref{DifEq1}) with $C^g_{\bf Y}(y)$ being replaced by $C^g_{\bf Y}(y)+ D(y)$. Similarly, one can show that the differential equation (\ref{DifEq1}) and Conditions \ref{con4}-\ref{con6} imply the integral equation (\ref{constar}) and Condition \ref{con41}.

To summarize, the current paper can be a starting point for the rigorous investigation of different types of the optimality equation for impulsively controlled PDMP.

\section*{Acknowledgements}

This research was supported by the Royal Society International Exchanges award IE160503, by FCT and CIDMA within project UID/MAT/04106/2013, and  by TOCCATA FCT project \linebreak PTDC/EEI-AUT/2933/2014.

\section{Appendix} \label{sec6}

\textit{Proof of Proposition \ref{prop0}.}

The function $V(\cdot)$ is lower semicontinuous by Theorem \ref{t1}. Then the function $G(\cdot)$ is lower semicontinuous, as seen in the proof of Theorem \ref{t1}. By Proposition 7.32 of \cite{Bertsekas:1978}, for each $x\in \textbf{X}$, $\inf_{a\in{\bf A}} G(x,\theta,a)$ defines a lower semicontinuous function on $\bar\RR^0_+$, and thus $\Theta(x)$ is closed and thus compact in $\bar\RR^0_+$. The nonemptyness of $\Theta(x)$ is by Theorem \ref{t1}.

By Proposition D.5 of \cite{las96}, $\inf_{a\in{\bf A}} G(x,\theta,a)$ defines a measurable function on $\textbf{X}\times \bar\RR^0_+.$  Then the graph of the multifunction $\Theta(\cdot)$, given by $\{(x,\theta)\in \textbf{X}\times \bar\RR^0_+: \inf_{a\in{\bf A}} G(x,\theta,a) =V(x)\}$, is measurable and hence the multifunction $\Theta(\cdot)$ is Borel-measurable by Proposition D.4 of \cite{las96}.
By Proposition D.5 of \cite{las96}, $\theta^\ast(x)=\inf_{\theta\in \Theta(x)}\theta$ defines a measurable function on $\textbf{X}.$  $~\hfill\Box$

{\it Proof of Proposition \ref{prop3}.}

According to Theorem \ref{t1} and inequalities $0\le {\cal V}(x)\le K$, it is sufficient to show only the uniqueness, namely, we will show that if $V$ is a bounded lower semicontinuous solution to (\ref{PZZeqn01}), then $V={\cal V}$.

Consider the obvious formula
\begin{eqnarray}
{\cal V}(x_0,\pi) &=& V(x_0)+\lim_{N\to\infty} \sum_{i=1}^N E^\pi_{x_0}\left[\vphantom{\int_X} \tilde C(X_{i-1},(\Theta_i,A_i))\right. \nonumber \\
&&\left.+\int_{\bf X} V(y) Q(dy|X_{i-1},(\Theta_i,A_i))-V(X_{i-1})\right]-\lim_{N\to\infty} E^\pi_{x_0}[V(X_N)]
\label{e5}
\end{eqnarray}
valid for each strategy $\pi$ such that
\begin{equation}\label{e6}
{\cal V}(x_0,\pi)\le K.
\end{equation}
Note that, for such strategies,
$$E^\pi_{x_0}[T_{stop}]=\sum_{i=1}^\infty P^\pi_{x_0}(T_{stop}\ge i)\le \frac{K}{\delta},$$
so that $\lim_{N\to\infty} P^\pi_{x_0}(T_{stop}\ge N)=0$
and
$$\lim_{N\to\infty} E^\pi_{x_0}[V(X_N)]\le \sup_{x\in{\bf X}} V(x)\lim_{N\to\infty} P^\pi_{x_0}(T_{stop}\ge N)=0.$$
Now the stationary deterministic strategy $(\varphi^*_{\theta},\varphi^*_a)$, providing the infimum in (\ref{PZZeqn01}), is uniformly optimal and
$${\cal V}^*(x_0)=\inf_{\pi} {\cal V}(x_0,\pi)={\cal V}(x_0,(\varphi^*_{\theta},\varphi^*_a))=V(x_0)$$
because all the other strategies except for those satisfying (\ref{e6}) cannot give smaller value for ${\cal V}(x_0,\pi)$.
\hfill $\Box$

{\it Proof of Proposition \ref{prop2}.}

It suffices to prove that for all $c>0$ the function $g(s) \defi h(s) +cs$ is nondecreasing on $[0,t]\cap\RR^0_+$.

Note that for any $s \in (0,\, t)$ there exists $\varepsilon =  \varepsilon_s$ such that
\vspace{1mm}

(i) either $g(s) < g(\sigma)$ for all $\sigma \in (s,\, s + \varepsilon)$,
\vspace{1mm}

(ii) or $g(\sigma) < g(s)$ for all $\sigma \in (s - \varepsilon,\, s)$.
\vspace{2mm}

Suppose that $g(s_1) > g(s_2)$ for some $0 \le s_1 < s_2 \le t$. Our aim is to come to a contradiction.

Take $y \in (g(s_2),\, g(s_1))$ and take $s_* = \inf A$, where $A = \{ s \in [s_1,\, s_2] : g(s) < y \}$. Note that $A$ contains $s_2$, and therefore is nonempty.

If $s_* = s_1$, then on each interval $[s_1,s_1+\delta]$ there are points from $A$, so that
$$
\underline{\lim}_{s \to s_1^+}\, g(s) \le y < g(s_1),
$$
in contradiction with the right lower semicontinuity of $g$ at $s_1$.

If $s_* = s_2$, we have
$$
g(s_2) < y \le \underline{\lim}_{s \to s_2^-}\, g(s) \le \overline{\lim}_{s \to s_2^-}\, g(s),
$$
in contradiction with the left upper semicontinuity of $g$ at $s_2$.

It follows that $s_1 < s_* < s_2$.
\vspace{1mm}

For all $s \in [s_1,\, s_*)$ we have $g(s) \ge y$, therefore
\begin{equation}\label{p2-1}
y \le \underline{\lim}_{s \to s_*^{-}}\, g(s) \le \overline{\lim}_{s \to s_*^{-}}\, g(s) \le g(s_*)
\end{equation}
because $g$ is left upper semicontinuous at $s_*$.
On the other hand, there exists a sequence $\{ s_k \} \subset A$ converging to $s_*$. If at least one term of this sequence coincides with $s_*$ then $g(s_*) < y$. If, otherwise, no terms of the sequence coincide with $s_*$ then all $s_k > s_*$ and, for all $k$, $g(s_k)<y$. Since $g$ is right lower semicontinuous at $s_*$,
\begin{equation}\label{p2-2}
g(s_*) \le \underline{\lim}_{s \to s_*^+}\, g(s) \le \underline{\lim}_{k \to \infty}\, g(s_k) \le y.
\end{equation}
It follows from (\ref{p2-1}) and (\ref{p2-2}) that $g(s_*) = y$, and both conditions (i), (ii) are violated for $s_*$: for all $s\in[s_1,s_*]$, $g(s)\ge g(s_*)$ and in each right neighbourhood of $s$ there are points $s\in A$ such that $g(s_*)>g(s)$. \hfill $\Box$

{\it Proof of Proposition \ref{prop1}.}

One easily sees that $h$ is continuous on $[0,\, t]$. Both $h$ and $-h$ satisfy the conditions of Proposition \ref{prop2}, therefore both $h$ and $-h$ are nondecreasing. It follows that $h$ is constant. \hfill $\Box$

{ {\it Proof of Proposition \ref{prqm}.}
Condition \ref{con2} follows from Condition \ref{as3}. According to Theorem \ref{th1}, the bounded lower semicontinuous non-negative function $V$ satisfies the Bellman equation (\ref{PZZeqn01}). The function ${\cal V}^*$ is bounded  because of Conditions \ref{asa}  and \ref{as3}. According to Remark \ref{rem9}, the Bellman equation cannot have another bounded lower semicontinuous solution. Therefore, $V={\cal V}^*$ is the minimal $\RR^0_+$-valued solution to equation (\ref{PZZeqn01}) and the strategy $(\varphi^*_\theta,\varphi^*_a)$ is uniformly optimal by Theorem \ref{t1}.
\hfill $\Box$}

\end{document}